\documentclass[12pt]{amsart}
  \usepackage{latexsym}
  \usepackage[all]{xy}
  \usepackage{amsthm}
 \numberwithin{equation}{section}
 \usepackage{color}
  \usepackage{hyperref}

  \newcommand{\Ra}{\Rightarrow}
\def\cK{\mathcal K}
\def\ot{{\otimes}}
\def\1c#1{\stackrel{#1}{\to}}
\def\2c#1{\stackrel{#1}{\Ra}}
\def\EM{\mathrm{EM}}
\def\Mnd{\mathrm{Mnd}}
\def\Cmd{\mathrm{Cmd}}
\def\ilift#1{\overrightarrow{#1}}
\def\plift#1{\overleftarrow{#1}}
\def\A{\EM(\cK)(I(l),J^w(s,\psi))}
\def\B{\EM(\cK)(I(l),{\widehat{st}})}
\def\nast{\! \ast \!}

  \newtheorem{proposition}{Proposition}[section]
  \newtheorem{lemma}[proposition]{Lemma}
  
  \newtheorem{corollary}[proposition]{Corollary}
  \newtheorem{theorem}[proposition]{Theorem}

  \theoremstyle{definition}
  \newtheorem{definition}[proposition]{Definition}
  \newtheorem{example}[proposition]{Example}
  \newtheorem{examples}[proposition]{Examples}

  \theoremstyle{remark}

  \newcounter{c}
  
  \newcommand{\etyk}[1]{\vspace{-7.4mm}$$\begin{equation}\Label{#1}
  \addtocounter{c}{1}}
  \renewcommand{\]}{\ifnum \value{c}=1 $$\else \end{equation}\fi}
\begin{document}

 \title{The weak theory of monads}
 \author{Gabriella B\"ohm}
\address{Research Institute for Particle and Nuclear Physics, Budapest,
\newline\indent H-1525
 Budapest 114, P.O.B.\ 49, Hungary}
\email{G.Bohm@rmki.kfki.hu}
\date{}
  \begin{abstract}
We construct a `weak' version $\EM^w(\cK)$ of Lack \&
Street's 2-category of monads in a 2-category $\cK$, by replacing
their compatibility constraint of 1-cells with the units of monads by an
additional condition on the 2-cells. 
A relation between monads in $\EM^w(\cK)$ and composite pre-monads in
$\cK$ is discussed. If $\cK$ admits Eilenberg-Moore constructions for monads, 
we define two symmetrical notions of {`weak liftings'} for monads in $\cK$.
If moreover idempotent 2-cells in $\cK$ split, we describe both kinds of
weak lifting via an appropriate pseudo-functor $\EM^w(\cK)\to \cK$. 
Weak entwining structures and partial entwining structures are shown to
realize weak liftings of a comonad for a monad in these respective senses. 
Weak bialgebras are characterized as algebras and coalgebras, such that the
corresponding monads weakly lift for the corresponding comonads and also the
comonads weakly lift for the monads.
\end{abstract}
  \maketitle

\section*{Introduction}

Many constructions, developed independently in Hopf algebra theory, turn out
to fit more general situations studied in category theory. 
For example, crossed products with a Hopf algebra in \cite{Swe}, \cite{BCM},
\cite{DT} are examples of a wreath product in \cite{LS}.
As another example, the comonad induced by the underlying coalgebra in a Hopf
algebra $H$, has a lifting to the category of modules over any $H$-comodule
algebra. So-called Hopf modules are comodules (also called coalgebras) for the
lifted comonad. Galois property of an algebra extension by a Hopf algebra
turns out to correspond to comonadicity of an appropriate functor \cite{Pepe}.  

Motivated by noncommutative differential geometry, these constructions were
extended from Hopf algebras to coalgebras (over a commutative base ring) and
to corings (over an arbitrary base ring), see e.g. the pioneering paper
\cite{Brz:coring}. The resulting theory turns out to fit the same categorical
framework, only the occurring (co)monads have slightly more complicated forms
\cite{Pepe}.  

For a study of non-Tannakian monoidal categories (i.e. those that admit no
strict monoidal fiber functor to the module category of some commutative ring),
another direction of generalization was proposed in \cite{BNSz}.
The essence of this approach is a weakening of the unitality of some maps 
and it leads to the replacing of a Hopf algebra by a `weak' Hopf algebra.
In the last decade many Hopf algebraic constructions were extended to the
weak setting. Weak crossed products were studied e.g. in \cite{CaeDeG}
\cite{WCP} and \cite{RodRap}. Weak Galois theory was developed, among other
papers, in \cite{CaeDeG}, \cite{Brz:coring}.
However, just because in these generalizations one deals with non-unital maps,
they do not fit the categorical framework of (co)monads, their wreath products
and liftings.
\bigskip

The aim of the current paper is to provide a categorical framework for 
`weak' constructions. 
For this purpose, in Section \ref{sec:EM^w} we construct, for any 2-category
$\cK$, a 2-category $\EM^w(\cK)$ which contains the 
2-category $\EM(\cK)$ in \cite{LS} as a vertically full 2-subcategory. 
In $\EM^w(\cK)$ 0-cells are the same as in $\EM(\cK)$, i.e. monads in $\cK$.
Since we aim to describe constructions in terms of non-unital maps, in the
definition of a 1-cell in $\EM^w(\cK)$ we impose the same compatibility
condition with the multiplications of monads which is required in $\EM(\cK)$,
but we relax the compatibility condition in $\EM(\cK)$ with the
units. Certainly, without compensating it with some other requirements, 
we would not obtain a 2-category. We show that imposing one further axiom on
the 2-cells in addition to the axiom in $\EM(\cK)$, $\EM^w(\cK)$ becomes a
2-category, with the same horizontal and vertical composition laws used in
$\EM(\cK)$.  

It was observed in \cite{CaeDeG} that smash products (and more generally
crossed products \cite{RodRap}) by weak bialgebras are not unital
algebras. Motivated by the definition of a pre-unit in \cite{CaeDeG}, we study
pre-monads (defined in Section \ref{sec:premonad}) in any 2-category $\cK$.
In Section \ref{sec:premonad} we interpret `weak crossed products' in
\cite{WCP} as monads in $\EM^w(\cK)$. 
This leads to a bijection between monads $t \1c s t$ in $\EM^w(\cK)$ and
pre-monad structures (in $\cK$) on the composite 1-cell $st$ with a `$t$-
linear' multiplication. 

Starting from Section \ref{sec:J^w}, we restrict our studies to 2-categories
$\cK$ which admit Eilenberg-Moore constructions for monads (in the sense of
\cite{Street}) and in which idempotent 2-cells split.
These assumptions are motivated by applications to bimodules. The bicategory
$\mathrm{BIM}$ of {\tt [Algebras; Bimodules; Bimodule maps]}, over a
commutative, associative and unital ring $k$, satisfies both
assumptions. However, in order to avoid technical complications caused by
non-strictness of the horizontal composition in a bicategory, we prefer to
restrict to 2-categories. In the examples, instead of the bicategory
$\mathrm{BIM}$, we can work with its image in the 2-category
$\mathrm{CAT}=${\tt [Categories; Functors; Natural transformations]}, under
the hom 2-functor $\mathrm{BIM}(k,-):\mathrm{BIM}\to \mathrm{CAT}$, which
image is a 2-category with the desired properties. 

For a 2-category $\cK$ which admits Eilenberg-Moore constructions for monads,
the inclusion 2-functor $I:\cK \to \EM(\cK)$ possesses a right
2-adjoint $J$, cf. \cite{LS}. In Section \ref{sec:J^w} we use the
splitting property of idempotent 2-cells in $\cK$ to construct a
factorization of $J$ through the inclusion 2-functor $\EM(\cK)\hookrightarrow
\EM^w(\cK)$ and an appropriate pseudo-functor $J^w:\EM^w(\cK)\to \cK$.
For a monad $t\1c s t$ in $\EM^w(\cK)$ and any 0-cell $k$ in $\cK$, we prove
that both monads $\cK(k,J^w(s))$ and $\cK(k, \widehat{st})$ in
$\mathrm{CAT}$ possess isomorphic Eilenberg-Moore categories, where
$\widehat{st}$ is a canonical retract monad of the pre-monad $st$.

In a 2-category $\cK$ which admits Eilenberg-Moore constructions for monads,
any monad $k \1c t k$ determines an adjunction $(k \1c f J(t), J(t) \1c v k)$
in $\cK$, cf. \cite{Street}. A lifting of a 1-cell $k\1c V k'$ in $\cK$ for
monads $k \1c t k$ and $k' \1c {t'} k'$ is, by definition, a 1-cell $J(t) \1c 
{\overline V} J(t')$, such that $v'{\overline V} = Vv$, cf. \cite{PowWat}.
In Section \ref{sec:lift} we define a `weak' lifting by replacing this
equality with the existence of a 2-cell $v'{\overline V} \2c\iota  Vv$,
possessing a retraction $Vv \2c\pi  v' {\overline V}$. 
This leads to two symmetrical notions of a weak lifting of a 2-cell $V \2c
\omega W$ for the monads $t$ and $t'$. A weak $\iota$-lifting ${\overline V}
\2c {\ilift \omega} {\overline W}$ is defined by the condition $\iota  \ast v' 
\ilift \omega = \omega v \ast\iota $ and a weak $\pi$-lifting ${\overline V}
\2c {\plift \omega} {\overline W}$ is defined by $v' \plift \omega \ast\pi
=\pi  \ast \omega v$. We show that any weak $\iota$-lifting and any weak
$\pi$-lifting of a 2-cell in $\cK$, if it exists, is isomorphic to the image
of an appropriate 2-cell in $\EM^w(\cK)$ under the pseudo-functor $J^w$. Both
a weak $\iota$-lifting and a weak $\pi$-lifting is proven to strictly 
preserve vertical composition and to preserve horizontal composition up to a
coherent isomorphism. We also give sufficient and necessary conditions for  
the existence of weak $\iota$-, and 
$\pi$-liftings of a 2-cell in $\cK$.   
\bigskip

A powerful tool to treat algebra extensions by weak bialgebras is provided by
`weak entwining structures' in \cite{CaeDeG}. A weak entwining structure in a
2-category $\cK$ consists of a monad $t$ and a comonad $c$, together with a
2-cell $tc \2c {} ct$ relating both structures in a way which generalizes a
mixed distributive law in \cite{Beck}. It was observed in \cite{Brz:coring}
that any weak entwining structure (in $\mathrm{BIM}$) induces a comonad
(called a `coring' in the particular case of the bicategory $\mathrm{BIM}$). 
In Section \ref{sec:appli} we show that -- in the same way as mixed
distributive laws in a 2-category $\cK$ provide examples of comonads in
$\EM(\cK)$ -- weak entwining structures provide examples of comonads in
$\EM^w(\cK)$. Moreover, if the 2-category $\cK$ satisfies the assumptions in
Section \ref{sec:J^w}, then the comonad in $\cK$, induced by a weak entwining 
structure, is an example of a weak $\iota$-lifting of a comonad for a monad. 

Studying partial coactions of Hopf algebras, in \cite{CaeJan:web} another
generalization of a mixed distributive law, a so called `partial entwining
structure' was introduced. Partial entwining structures (in $\mathrm{BIM}$)
were proven to induce comonads as well.  
We show that also partial entwining structures in a 2-category $\cK$ provide
examples of a comonad in $\EM^w(\cK)$. Moreover, if the 2-category
$\cK$ satisfies the assumptions in Section \ref{sec:J^w}, then the comonad in
$\cK$, induced by a partial entwining structure, is an example of a weak
$\pi$-lifting of a comonad for a monad. 

As a final application, weak bialgebras are characterized via weak liftings.
If a module $H$, over a commutative, associative and unital ring $k$, possesses
both an algebra and a coalgebra structure, then it induces two monads $t_R =
(-)\otimes_k H$ and $t_L=H \otimes_k (-)$, and two comonads $c_R = (-)
\otimes_k H$ and $c_L = H \otimes_k (-)$ on the category of $k$-modules. 
We relate weak bialgebra structures of $H$ to 
weak $\iota$-liftings of $c_R$ and $c_L$ for $t_R$ and $t_L$, respectively, and 
weak $\pi$-liftings $t_R$ and $t_L$ for $c_R$ and $c_L$, respectively.

{\bf Notations.}
We assume that the reader is familiar with the theory of 2-categories. For a
review of the occurring notions (such as a 2-category, a 2-functor and a
2-adjunction, monads, adjunctions and Eilenberg-Moore construction in a
2-category) we refer to the article \cite{KeSt}.

In a 2-category $\cK$, horizontal composition is denoted by juxtaposition and
vertical composition is denoted by $\ast$, 1-cells are represented by an arrow
$\1c{}$ and 2-cells are represented by $\2c {}$. 

For any 2-category $\cK$, $\Mnd(\cK)$ denotes the 2-category of monads in
$\cK$ as in \cite{Street} and $\Cmd(\cK):=\Mnd(\cK_*)_*$ denotes the
2-category of comonads in $\cK$, where $(-)_*$ refers to the vertical opposite
of a 2-category.
We denote by $\EM(\cK)$ the extended 2-category of monads in \cite{LS}. 
We use the reduced form of 2-cells in $\EM(\cK)$, see \cite{LS}.

\section{The 2-category $\EM^w(\cK)$}\label{sec:EM^w}

For any 2-category $\cK$, the 2-category $\EM^w({\cK})$ introduced in this
section extends the 2-category $\EM({\cK})$ in \cite{LS}.  

\begin{theorem} \label{thm:EMwK}
For any 2-category $\cK$, the following data constitute a 2-category,
to be denoted by $\EM^w(\cK)$.

\begin{itemize}
\item[] {\underline{$0$-cells}} are monads 
$(k \1c t k, tt \2c \mu t, k \2c \eta  t)$ in $\cK$.
\item[] {\underline{$1$-cells}} $(t,\mu,\eta )\to (t',\mu',\eta ')$ are pairs
  $(V,\psi)$, consisting of a 1-cell $k \1c V k'$ and a 2-cell $t'V \2c \psi V
  t$ in $\cK$, such that  
\begin{equation}\label{eq:1-cell}
V\mu \ast \psi t \ast t' \psi = \psi \ast \mu' V.
\end{equation}
The identity 1-cell is $t \stackrel{(k,t)}{\longrightarrow} t$.
\item[] {\underline{$2$-cells}} $(V,\psi) \Rightarrow (W,\phi)$ are 2-cells
  $V \2c \varrho Wt$ in $\cK$, such that 
\begin{eqnarray}
&&W\mu \ast \varrho t \ast \psi = W\mu \ast \phi t \ast t'\varrho
\label{eq:2-cell_a} \\
&&\varrho = W\mu \ast \phi t \ast \eta ' W t \ast \varrho.
\label{eq:2-cell_b}
\end{eqnarray}
The identity 2-cell is $(W,\phi)\2c {\phi \ast \eta ' W}
(W,\phi)$. 
\item[] {\underline{Horizontal composition}} of 2-cells 
$(V,\psi)  \2c \varrho (W,\phi)$,
$(V',\psi')  \2c \varrho' (W',\phi')$
(for 1-cells $(V,\psi), (W,\phi):t \to t'$ and $(V',\psi'), (W',\phi'):t' \to
  t''$) is given by 
\begin{equation}\label{eq:hor}
\varrho'\circ \varrho:=W'W\mu \ast W'\varrho t \ast W'\psi \ast \varrho' V.
\end{equation}
\item[] {\underline{Vertical composition}} of 2-cells
  $(V,\psi)\2c \varrho (W,\phi)\2c \tau (U,\Theta)$ 
(for 1-cells $(V,\psi), (W,\phi)$ and $(U,\Theta):t \1c{} t'$)
is given by 
\begin{equation}\label{eq:vert}
\tau \bullet \varrho:=U\mu \ast \tau t \ast \varrho.
\end{equation}
\end{itemize}
\end{theorem}

\begin{proof}
We verify only those axioms whose proof is different from the proof of the
respective axiom for $\EM(\cK)$.

The vertical composite of 2-cells $(V,\psi) \2c \varrho (W,\phi)\2c \tau
(U,\Theta)$ in $\EM^w(\cK)$ is checked to satisfy \eqref{eq:2-cell_a} by the
same computation used to verify that the vertical composite of 2-cells in
$\EM(\cK)$ is a 2-cell in $\EM(\cK)$. In order to see that $\tau \bullet
\varrho$ satisfies \eqref{eq:2-cell_b}, use the interchange law in $\cK$,
associativity of the multiplication $\mu$ of the monad $t$ and the fact that
$\tau$ satisfies \eqref{eq:2-cell_b}:
\begin{eqnarray*}
U\mu \ast \Theta t\ast \eta ' U t \ast (\tau \bullet \varrho) &=&
U\mu \ast \Theta t\ast \eta ' U t \ast U\mu \ast \tau t \ast \varrho \\
&=&U\mu \ast U\mu t \ast \Theta tt \ast \eta ' U tt \ast \tau t \ast \varrho = 
U\mu  \ast \tau t \ast \varrho = \tau \bullet \varrho.
\end{eqnarray*}
Associativity of the vertical composition follows by the same reasoning as in
the case of $\EM(\cK)$. 
Using the constraint \eqref{eq:1-cell}, it follows for any 1-cell $t
\stackrel{(W,\phi)}{\to} t'$ in $\EM^w(\cK)$ that
\begin{eqnarray*}
&&W\mu  \ast \phi t \ast t' \phi \ast t' \eta ' W =
\phi \ast \mu ' W \ast t' \eta ' W =
\phi
\quad \textrm{and}\\
&&W\mu  \ast \phi t \ast \eta ' Wt \ast \phi =
W\mu  \ast \phi t \ast t'\phi \ast \eta ' t' W =
\phi \ast \mu ' W \ast \eta ' t' W = \phi.
\end{eqnarray*}
Thus the 2-cell
$(W,\phi)\2c {\phi\ast \eta 'W} (W,\phi)$ in $\cK$ satisfies
\eqref{eq:2-cell_a} and \eqref{eq:2-cell_b}, proving that it is a 2-cell in
$\EM^w(\cK)$. 
It is immediate by condition \eqref{eq:2-cell_b} that for any 2-cell $(V,\psi)
\2c \varrho (W,\phi)$ in $\EM^w(\cK)$,
$$
(\phi \ast \eta ' W)\bullet \varrho =
W\mu  \ast \phi t \ast \eta ' W t \ast \varrho = \varrho.
$$
Using \eqref{eq:2-cell_a}, the interchange law in $\cK$ and then
\eqref{eq:2-cell_b}, one checks that also 
$$
\varrho \bullet (\psi \ast \eta ' V) =
W\mu  \ast \varrho t \ast \psi \ast \eta ' V =
W\mu  \ast \phi t \ast t' \varrho \ast \eta ' V =
W\mu  \ast \phi t \ast \eta ' Wt \ast \varrho =
\varrho.
$$
Hence there are identity 2-cells of the stated form.

On identity 2-cells 
$(V,\psi)\2c {\psi \ast \eta ' V} (V,\psi)$ and 
$(V',\psi') \2c {\psi' \ast \eta '' V'} (V',\psi')$, 
the horizontal composite comes out as
$$
\left(\psi'\ast \eta '' V'\right) \circ (\psi \ast \eta ' V)
=V' \psi \ast \psi' V \ast \eta '' V' V,
$$
where we applied \eqref{eq:1-cell} for $\psi$ and unitality of the monad $t'$.
That is to say, the horizontal composite of the 1-cells $t
\stackrel{(V,\psi)}{\longrightarrow} t' \stackrel{(V',\psi')}{\longrightarrow}
t''$ is the 1-cell $(V'V, V' \psi \ast \psi' V)$.
By the same computations used in the case of $\EM(\cK)$, the horizontal
composite of 2-cells $(V,\psi)  \2c \varrho (W,\phi)$,
$(V',\psi')  \2c \varrho' (W',\phi')$
is checked to satisfy condition \eqref{eq:2-cell_a}. It satisfies also
\eqref{eq:2-cell_b}, as 
\begin{eqnarray*}
W'W\mu  &\ast& W'\phi t \ast \phi' W t \ast \eta '' W'Wt\ast (\varrho' \circ
\varrho)\\
&=& W'W\mu  \ast W'\phi t \ast \phi' W t \ast \eta '' W'Wt\ast W'W\mu  \ast W'
\varrho t \ast W' \psi \ast \varrho' V \\
&=& W'W\mu  \ast W'W\mu t \ast W'\phi tt \ast W' t' \varrho t \ast \phi' V t
\ast \eta '' W' V t \ast  W' \psi \ast \varrho' V \\
&=& W'W\mu  \ast W'W\mu t \ast W' \varrho tt \ast W' \psi t \ast \phi' Vt \ast
\eta '' W' V t \ast W' \psi \ast \varrho' V \\
&=& W'W\mu  \ast W' \varrho t \ast W' V \mu  \ast W' \psi t \ast W' t' \psi
\ast \phi' t' V \ast \eta '' W' t' V \ast \varrho' V \\
&=& W'W\mu  \ast W' \varrho t \ast W'\psi \ast W' \mu ' V \ast \phi' t' V \ast
\eta '' W' t' V \ast \varrho' V \\
&=& W'W\mu  \ast W' \varrho t \ast W'\psi \ast  \varrho' V = \varrho' \circ
\varrho. 
\end{eqnarray*}
The first and last equalities follow by \eqref{eq:hor}.
In the second and fourth equalities we used the interchange law in $\cK$ and
associativity of the multiplication $\mu $ of the monad $t$.
In the third equality we used that $\varrho$ satisfies \eqref{eq:2-cell_a}. 
The fifth equality is derived by using that $\psi$ satisfies
\eqref{eq:1-cell}. The penultimate equality follows by using that $\varrho'$
satisfies \eqref{eq:2-cell_b}. 
This proves that the horizontal composite of 2-cells is a 2-cell.
Associativity of the horizontal composition in $\EM^w(\cK)$ is checked in the
same way as it is done in $\EM(\cK)$. Obviously, for any 2-cell
$(V,\psi) \2c \varrho (W,\phi)$, 
$\varrho \circ \eta  =
W\mu  \ast W \eta  t \ast \varrho =
\varrho$.
By \eqref{eq:2-cell_a} and \eqref{eq:2-cell_b}, also
$\eta '\circ \varrho 
= W\mu  \ast \varrho t \ast \psi \ast \eta ' V 
=\varrho$.
Hence the identity 2-cell $(k,t) \2c \eta  (k,t)$ is a unit for the horizontal
composition, proving the stated form $(k,t)$ of the identity 1-cell $t\to t$.

The interchange law in $\EM^w(\cK)$ is checked in the same way as it is done in
$\EM(\cK)$. 
\end{proof}

Clearly, any 1-cell in $\EM(\cK)$ is a 1-cell also in $\EM^w(\cK)$. 
For a 1-cell $t \1c {(W,\phi)} t'$ in $\EM(\cK)$, any 2-cell $\varrho$ in
$\cK$ of target $Wt$ satisfies \eqref{eq:2-cell_b}. Hence 2-cells in
$\EM^w(\cK)$ between 1-cells of $\EM(\cK)$ are the same as the 2-cells in
$\EM(\cK)$. Comparing the formulae of the horizontal and vertical compositions
in $\EM(\cK)$ and $\EM^w(\cK)$, we conclude that $\EM(\cK)$ is a vertically
full 2-subcategory of $\EM^w(\cK)$.

One may ask what 2-subcategory of $\EM^w(\cK)$ plays the role of the
2-subcategory, obtained as an image of $\Mnd(\cK)$ in $\EM(\cK)$. As the lemmata
below show, there seems to be no unique answer to this question. This is
because, for some 1-cells $t\1c {(V,\psi)} t'$ and $t \1c {(W,\phi)} t'$ in
$\EM^w(\cK)$ and a 2-cell $V \2c \omega W$ in $\cK$, the 2-cells 
$\omega t \ast \psi \ast \eta 'V$ and $\phi\ast \eta ' W \ast \omega$ 
in $\cK$ need not be equal. Still, there are two distinguished sets
of 2-cells in $\EM^w(\cK)$ on equal footing, both closed under the horizontal
and vertical compositions and both containing the identity 2-cells.

\begin{lemma}\label{lem:2-cells}
For any 2-category $\cK$,
let $t\1c {(V,\psi)} t'$ and $t \1c {(W,\phi)} t'$ be 1-cells in $\EM^w(\cK)$
and $V \2c \omega W$ be a 2-cell in $\cK$. 

\begin{itemize}
\item[{(1)}] The following assertions are equivalent.
\begin{itemize}
\item[{(i)}] $\omega t \ast \psi \ast \eta ' V:(V,\psi) \2c {}
  (W,\phi)$ is a 2-cell in $\EM^w(\cK)$.
\item[{(ii)}] $\omega t\ast \psi = W\mu  \ast \phi t \ast
  t'\omega  t \ast t' \psi \ast t' \eta ' V$;
\item[{(iii)}] $W\mu  \ast \phi t \ast \eta ' Wt \ast \omega t
  \ast \psi \ast \eta ' V =\omega t \ast \psi \ast \eta ' V$  
and \\
$W\mu  \ast \phi t \ast \eta ' Wt \ast \omega t \ast \psi   
=W \mu  \ast \phi t \ast t' \omega t \ast t' \psi \ast t' \eta ' V$;
\end{itemize}

\item[{(2)}] The following assertions are equivalent.
\begin{itemize}
\item[{(i)}] $\phi\ast \eta ' W \ast \omega:(V,\psi) \2c {}
  (W,\phi)$ is a 2-cell in $\EM^w(\cK)$. 
\item[{(ii)}] $\phi \ast t' \omega = W\mu  \ast \phi t \ast
  \eta ' Wt \ast \omega t \ast\psi$;  
\item[{(iii)}] $W\mu  \ast \phi t \ast \eta ' Wt \ast \omega t
  \ast \psi \ast \eta 'V = \phi \ast \eta ' W \ast \omega$   
and \\
$W\mu  \ast \phi t \ast \eta ' Wt \ast \omega t \ast \psi   
=W \mu  \ast \phi t \ast t' \omega t \ast t' \psi \ast t' \eta ' V$;
\end{itemize}

\item[{(3)}] The following assertions are equivalent.
\begin{itemize}
\item[{(i)}] $\phi \ast t' \omega = \omega t \ast \psi$; 
\item[{(ii)}] $\phi\ast \eta ' W \ast \omega$ and $\omega t \ast \psi \ast
  \eta ' V$ are (necessarily equal) 2-cells $(V,\psi) \2c {} (W,\phi)$ in
  $\EM^w(\cK)$. 
\end{itemize}
\end{itemize}
\end{lemma}

\begin{proof}
\underline{(1)} \underline{(i)$\Leftrightarrow$ (ii)} 
Using that $\psi$ satisfies \eqref{eq:1-cell} together with the unitality of
the monad $t'$, condition \eqref{eq:2-cell_a} for $\varrho:= \omega t \ast
\psi \ast \eta ' V$ comes out as the equality in part (ii).
Hence in order to prove the equivalence of assertions (i) and
(ii), we need to show that (ii) implies that
$\varrho$ satisfies \eqref{eq:2-cell_b}. Indeed, applying the equality in part
(ii) in the second step, we obtain 
\begin{equation} \label{eq:omega_b}
W\mu  \ast \phi t \ast \eta ' Wt \ast \omega t \ast \psi \ast \eta ' V
=W\mu  \ast \phi t \ast t' \omega t \ast t' \psi \ast t' \eta ' V \ast \eta ' V
= \omega t \ast \psi \ast \eta ' V.
\end{equation}

\underline{(ii)$\Leftrightarrow$(iii)}
We have seen in the proof of equivalence (i)$\Leftrightarrow$(ii) above, 
that assertion (ii) implies \eqref{eq:omega_b}, i.e. the first
condition in part (iii). The second condition is checked as
follows. 
\begin{eqnarray}
W\mu  \ast \phi t \ast \eta ' Wt \ast \omega t \ast \psi 
&=& W\mu  \ast \phi t \ast \eta ' Wt \ast \omega t \ast V\mu  \ast \psi t \ast
t' \psi \ast \eta ' t' V \nonumber\\
&=& W\mu  \ast W\mu t \ast \phi tt \ast \eta ' W tt \ast \omega tt \ast \psi t
\ast \eta ' Vt \ast \psi\nonumber\\
&=& W\mu  \ast \omega tt \ast \psi t \ast \eta ' Vt \ast \psi\nonumber\\
&=& \omega t \ast \psi
= W \mu  \ast \phi t \ast t' \omega t \ast t' \psi \ast t' \eta ' V.
\label{eq:omega_a}
\end{eqnarray}
The first equality follows by \eqref{eq:1-cell} and unitality of the monad
$t'$. 
In the second equality we used associativity of $\mu $ and the interchange law.
The third equality is obtained by using \eqref{eq:omega_b}.
The penultimate equality follows by \eqref{eq:1-cell} and unitality of the
monad $t'$ again. 
The last equality follows by assertion (ii).

Conversely, assume that the identities in part (iii) hold. Then
\begin{eqnarray*}
\omega t \ast \psi
&=& \omega t \ast V\mu  \ast \psi t \ast t' \psi \ast \eta ' t' V
= W\mu  \ast W\mu t \ast \phi tt \ast \eta ' Wtt \ast \omega tt \ast \psi t
\ast \eta ' Vt \ast \psi\\
&=& W\mu  \ast \phi t \ast \eta ' Wt \ast \omega t \ast \psi
= W\mu  \ast \phi t \ast t' \omega t \ast t'\psi \ast t' \eta ' V.
\end{eqnarray*}
The first equality follows by \eqref{eq:1-cell} and unitality of the monad
$t'$. 
In the second equality we applied the first condition in part (iii). The third
equality is obtained by using the associativity of $\mu $, 
the interchange law and \eqref{eq:1-cell} again.
The last equality follows by the second condition in part (iii).

\underline{(2)} \underline{(i)$\Leftrightarrow$(ii)} 
Using that $\phi$ satisfies \eqref{eq:1-cell} together with the unitality of
the monad $t'$, condition \eqref{eq:2-cell_a} for $\varrho:= \phi \ast \eta ' W
\ast \omega$ comes out as the equality in part (ii). 
Condition \eqref{eq:2-cell_b} holds for $\varrho$ automatically,
i.e. it follows by applying \eqref{eq:1-cell} for $\phi$.

\underline{(ii)$\Leftrightarrow$(iii)}
If (iii) holds, then 
\begin{eqnarray*}
\phi \ast t' \omega 
&=& W\mu  \ast \phi t \ast t' \phi \ast t' \eta ' W \ast t' \omega\\
&=& W\mu  \ast \phi t \ast t' W\mu  \ast t' \phi t \ast t' \eta ' Wt \ast t'
\omega t \ast t' \psi \ast t' \eta ' V\\
&=& W\mu  \ast \phi t \ast t' \omega t \ast t'\psi \ast t' \eta ' V
=W\mu  \ast \phi t \ast \eta ' Wt \ast \omega t \ast \psi.
\end{eqnarray*}
The first equality follows by applying \eqref{eq:1-cell} for $\phi$, together
with the unitality of the monad $t'$.
The second equality is obtained by the first identity in part (iii). 
In the penultimate equality we applied the interchange law, associativity of
$\mu $, \eqref{eq:1-cell} on $\phi$ and unitality of $t'$.
The last equality follows by the second condition in part (iii).

Conversely, if assertion (ii) holds, then the first condition in part (iii) is
proven by composing both sides of the equality in part (ii) by $\eta ' V$ on
the right.  
The second condition, i.e. equality \eqref{eq:omega_a}, is proven by the
following computation.  
\begin{eqnarray*}
W\mu  \ast \phi t \ast t' \omega t \ast t' \psi \ast t' \eta ' V
&=& W\mu  \ast W\mu t \ast \phi tt \ast t' \phi t \ast t' \eta ' Wt \ast t'
\omega t \ast t' \psi \ast t' \eta ' V\\
&=& W\mu  \ast \phi t \ast t'\phi \ast t' \eta ' W \ast t' \omega \\
&=& \phi \ast t' \omega
= W\mu  \ast \phi t \ast \eta ' Wt \ast \omega t \ast \psi.
\end{eqnarray*}
The first equality follows by applying \eqref{eq:1-cell} for $\phi$, together
with the unitality of the monad $t'$.
The second equality is obtained by using the associativity of $\mu $, the
interchange law and the first condition in part (iii).
The third equality follows by applying \eqref{eq:1-cell} for $\phi$, together
with the unitality of the monad $t'$ again. The last equality follows by part
(ii). 

\underline{(3)}
Assume first that assertion (3)(ii) holds. Then
by parts (1) and (2), also (1)(ii) and (2)(ii)
hold, implying 
\eqref{eq:omega_a}. Hence
$$
\omega t \ast \psi 
= W\mu  \ast \phi t \ast t' \omega t \ast t' \psi \ast t' \eta ' V' 
= W\mu  \ast \phi t \ast \eta ' Wt \ast \omega t \ast \psi 
= \phi \ast t' \omega.
$$

Conversely, if assertion (3)(i) holds, then 
\begin{eqnarray*}
&&W\mu  \ast \phi t \ast t' \omega t \ast t' \psi \ast t' \eta ' V
= W\mu  \ast \omega tt \ast \psi t \ast t' \psi \ast t' \eta ' V
= \omega t \ast \psi, \quad \textrm{and}\\
&&W\mu  \ast \phi t \ast \eta ' Wt \ast \omega t \ast \psi
= W\mu  \ast \phi t \ast \eta ' Wt \ast \phi \ast t' \omega 
= \phi \ast t' \omega,
\end{eqnarray*}
where in both computations the first equality follows by (3)(i) and the second
equality follows by \eqref{eq:1-cell} and the unitality of the monad $t'$.
We conclude by parts (1) and (2) that both $\omega t \ast \psi \ast \eta ' V$
and $\phi \ast \eta ' W \ast \omega$ are 2-cells in $\EM^w(\cK)$.
It follows by comparing the first identities in (1)(iii) and (2)(iii) that  
$\phi \ast \eta ' W \ast \omega = 
\omega t \ast \psi \ast \eta ' V$, as stated.
\end{proof}

Next we investigate the behaviour of the correspondences in Lemma
\ref{lem:2-cells} with respect to the horizontal and vertical compositions in
$\cK$ and $\EM^w(\cK)$.

\begin{lemma}\label{lem:lift_comp}
For any 2-category $\cK$, 
let $(V,\psi)$, $(W,\phi)$ and $(U,\theta)$ be 1-cells $t \1c {} t'$ and $(V',
\psi')$ and $(W', \phi')$ be 1-cells $t' \1c {} t''$ in $\EM^w(\cK)$. 
\begin{itemize}
\item[{(1)}] If some 2-cells $V\2c \omega W$ and $V' \2c {\omega '} W'$ in
  $\cK$ satisfy the equivalent conditions in Lemma \ref{lem:2-cells}(1),
  then 
$\left(\omega' t' \ast \psi' \ast \eta '' V'\right) \circ 
  \left(\omega t \ast \psi  \ast \eta ' V \right) = 
\omega' \omega t \ast V' \psi \ast \psi' V \ast \eta ''V'V$. Hence in particular
also $\omega' \omega$ satisfies the equivalent conditions in Lemma
\ref{lem:2-cells}(1).
\item[{(2)}] If some 2-cells $V\2c \omega W$ and $V' \2c {\omega '} W'$ in
  $\cK$ satisfy the equivalent conditions in Lemma \ref{lem:2-cells}(2),
  then 
$\left(\phi'\ast \eta '' W' \ast \omega'\right)\circ 
 \left(\phi \ast \eta ' W  \ast \omega \right)
=W'\phi\ast \phi' W \ast \eta '' W' W \ast \omega ' \omega$. Hence in particular
also $\omega' \omega$ satisfies the equivalent conditions in Lemma
\ref{lem:2-cells}(2).
\item[{(3)}] If some 2-cells $V\2c \omega W \2c {\kappa} U$ in
  $\cK$ satisfy the equivalent conditions in Lemma \ref{lem:2-cells}(1),
  then 
$\left(\kappa t \ast \phi \ast \eta ' W\right) \bullet 
  \left(\omega t \ast \psi  \ast \eta ' V \right) = 
\kappa t \ast \omega t \ast \psi \ast \eta ' V$. Hence in particular
also $\kappa\ast \omega$ satisfies the equivalent conditions in Lemma
\ref{lem:2-cells}(1).
\item[{(4)}] If some 2-cells $V\2c \omega W \2c {\kappa} U$ in
  $\cK$ satisfy the equivalent conditions in Lemma \ref{lem:2-cells}(2),
  then 
$\left(\theta \ast \eta ' U \ast \kappa\right)\bullet 
 \left(\phi \ast \eta ' W  \ast \omega \right)=
\theta \ast \eta ' U  \ast \kappa \ast \omega$. Hence in particular
also $\kappa\ast \omega$ satisfies the equivalent conditions in Lemma
\ref{lem:2-cells}(2).
\end{itemize}
\end{lemma}

\begin{proof}
\underline{(1)} This compatibility with the horizontal composition follows by
applying \eqref{eq:1-cell} for $\psi$, and unitality of the monad $t'$.

\underline{(2)} This follows by using that $\omega$ obeys Lemma
\ref{lem:2-cells}(2)(ii). 

\underline{(3)} This compatibility with the vertical composition follows 
using that, by Lemma \ref{lem:2-cells}(1), $\omega t \ast \psi \ast \eta ' V$
satisfies \eqref{eq:2-cell_b}.  

\underline{(4)} This  assertion follows by using that $\kappa$ satisfies Lemma
\ref{lem:2-cells}(2)(ii). 
\end{proof}

The message of Lemma \ref{lem:2-cells} and Lemma \ref{lem:lift_comp} can be
summarized as follows.

\begin{corollary}\label{cor:Mnd}
Consider an arbitrary 2-category $\cK$.
\begin{itemize}
\item[{(1)}] There is a 2-category, to be denoted by $\Mnd^\iota (\cK)$,
  defined by  the following data. 
\begin{itemize}
\item[] \underline{0-cells} are monads $t$ in $\cK$;
\item[] \underline{1-cells} $t\1c{(V,\psi)}t'$ are the same as 1-cells in
$\EM^w(\cK)$, cf. \eqref{eq:1-cell};
\item[] \underline{2-cells} $(V,\psi)\2c{\omega} (W,\phi)$ are 2-cells
$V\2c{\omega}W$ in $\cK$, satisfying the equivalent conditions in Lemma
\ref{lem:2-cells}(1);
\item[] \underline{horizontal and vertical compositions} are the same as in
  $\cK$. 
\end{itemize}
Furthermore, there is a 2-functor $G^\iota :\Mnd^\iota (\cK)\to \EM^w(\cK)$,
acting on the 0-, and 1-cells as the identity map and taking a 2-cell
$(V,\psi)\2c{\omega} (W,\phi)$ to $\omega t \ast \psi \ast \eta ' V$. 

\item[{(2)}] There is a 2-category, to be denoted by $\Mnd^\pi (\cK)$, defined
  by the following data.  
\begin{itemize}
\item[] \underline{0-cells} are monads $t$ in $\cK$; 
\item[] \underline{1-cells} $t\1c{(V,\psi)}t'$ are the same as 1-cells in
$\EM^w(\cK)$, cf. \eqref{eq:1-cell};
\item[] \underline{2-cells} $(V,\psi)\2c{\omega} (W,\phi)$ are 2-cells
$V\2c{\omega}W$ in $\cK$, satisfying the equivalent conditions in Lemma
\ref{lem:2-cells}(2);
\item[] \underline{horizontal and vertical compositions} are the same as in
  $\cK$. 
\end{itemize}
Furthermore, there is a 2-functor $G^\pi:\Mnd^\pi (\cK)\to \EM^w(\cK)$, acting
on the 0-, and 1-cells as the identity map and taking a 2-cell
$(V,\psi)\2c{\omega} (W,\phi)$ to $\phi\ast \eta ' W \ast \omega$.  
\end{itemize}
\end{corollary}

Clearly, both $\Mnd^\iota (\cK)$ and $\Mnd^\pi (\cK)$ contain $\Mnd(\cK)$ as a
vertically full subcategory.

\section{Monads in $\EM^w(\cK)$ and pre-monads in $\cK$}\label{sec:premonad} 

It was observed in \cite[page 257]{LS} that monads in $\EM(\cK)$ induce
monads in $\cK$. The aim of this section is to give a
similar interpretation of monads in $\EM^w(\cK)$.

A monad in $\EM^w(\cK)$ is given by a triple $((s,\psi),\nu,\vartheta)$,
consisting of a 1-cell $(t,\mu ,\eta ) \1c {(s,\psi)} (t,\mu ,\eta )$ and
2-cells $(s,\psi)\circ (s,\psi) \2c \nu (s,\psi)$ and $(k,t) \2c \vartheta
(s,\psi)$ in $\EM^w(\cK)$, such that  
\begin{eqnarray*}
\nu \bullet (\nu \circ (s,\psi))&=& 
\nu \bullet ((s,\psi)\circ \nu)\\
\nu \bullet (\vartheta\circ (s,\psi))=&(s,\psi)&=
\nu \bullet ((s,\psi)\circ \vartheta).
\end{eqnarray*}
In light of Theorem \ref{thm:EMwK}, this means a 1-cell $k \1c s k$, and
2-cells $ts \2c \psi st$, $ss \2c \nu st$ and $k \2c \vartheta st$ in $\cK$,
subject to the following identities.
\begin{eqnarray}
&&\psi \ast \mu s = s\mu  \ast \psi t \ast t \psi,
\label{eq:psi_1}\\
&&s\mu  \ast \psi t \ast t \nu = s\mu  \ast \nu t \ast s \psi \ast \psi s,
\label{eq:mu_2_a}\\
&&s\mu  \ast \psi t \ast \eta st \ast \nu =\nu,
\label{eq:mu_2_b}\\
&&s\mu  \ast \psi t \ast t \vartheta = s\mu  \ast \vartheta t,
\label{eq:eta_2_a}\\
&&s\mu  \ast \nu t \ast s\nu = s\mu  \ast \nu t \ast s \psi \ast \nu s,
\label{eq:mu_associ}\\
&&s\mu  \ast \nu t \ast s \psi \ast \vartheta s = \psi \ast \eta s,
\label{eq:eta_lun}\\
&&s\mu  \ast \nu t \ast s \vartheta = \psi \ast \eta s.
\label{eq:eta_run}
\end{eqnarray}
Condition \eqref{eq:psi_1} expresses the requirement that $(s,\psi)$ is a
1-cell in $\EM^w(\cK)$,
\eqref{eq:mu_2_a} and \eqref{eq:mu_2_b} together mean that $\nu$ is a 2-cell
in $\EM^w(\cK)$, 
\eqref{eq:eta_2_a} means that $\vartheta$ is a 2-cell in $\EM^w(\cK)$ (condition
\eqref{eq:2-cell_b} on $\vartheta$ follows by the interchange law in $\cK$,
\eqref{eq:eta_2_a} and unitality of the monad $t$). Conditions
\eqref{eq:mu_associ}, \eqref{eq:eta_lun} and 
\eqref{eq:eta_run} express associativity and unitality of the monad
$((s,\psi),\nu,\vartheta)$, after being simplified using \eqref{eq:psi_1},
\eqref{eq:mu_2_a}, \eqref{eq:mu_2_b} and \eqref{eq:eta_2_a}.   

Note that a monad $(t\1c{(s,\psi)}t,\nu,\vartheta)$ in $\EM^w(\cK)$, for a monad
$k\1c t k$ in $\cK$, is identical to a
`crossed product system' $(t,s,\psi,\nu)$ in the monoidal category
$\cK(k,k)$, in the sense of \cite[Definition 3.5]{WCP}, subject to the
`twisted' and `cocycle' conditions in \cite[Definition 3.3 and Definition
  3.6]{WCP}, the normalization condition in \cite[Proposition 3.7]{WCP} and
identities (11), (12) and (13) in \cite[Theorem 3.11]{WCP}, for $\vartheta$.

The following definition is inspired by \cite[Section 3.1]{CaeDeG} (see also
\cite[Definition 2.3]{WCP}).

\begin{definition}\label{def:pre_monad}
A {\em pre-monad} in a 2-category $\cK$ is a triple $(t,\mu ,\eta )$,
consisting of a 1-cell $k \1c t k$ and 2-cells $tt \2c \mu  t$ and $k \2c \eta
t$, such that the following conditions hold.
\begin{eqnarray}
&& \mu  \ast \mu t = \mu  \ast t\mu  \label{eq:pre_associ}\\
&& \mu  \ast \eta t = \mu  \ast t\eta  \label{eq:preun_1}\\
&& \mu  \ast \eta \eta  = \eta  \label{eq:preun_2}\\
&& \mu  \ast \mu t \ast \eta  tt =\mu . \label{eq:norm}
\end{eqnarray}
\end{definition}

Note that if $k\1c t k$ is a 1-cell in a 2-category $\cK$ and some 2-cells $tt
\2c \mu  t$ and $k \2c \eta  t$ satisfy $\mu \ast \mu t = \mu  \ast t\mu $ and
$\mu  \ast \eta t = \mu \ast t\eta  = \mu \ast \mu t \ast \eta \eta t$ as in
\cite[Section 3.1]{CaeDeG}, then $(t, \mu ':= \mu \ast \mu t \ast \eta tt,
\eta ':= \mu \ast \eta \eta )$ is a pre-monad in the sense of Definition 
\ref{def:pre_monad}. 

The motivation for a study of pre-monads stems from the following observation.

\begin{lemma}\label{lem:pre_monad}
Consider a pre-monad $(t,\mu ,\eta )$ in an arbitrary 2-category $\cK$.
\begin{itemize}
\item[{(1)}] 
  The 2-cell $\mu \ast \eta t$ is idempotent.
\item[{(2)}] 
  Assume that there exists a 1-cell ${\widehat t}$ and 2-cells $t \2c\pi 
  {\widehat t}$ and ${\widehat t}\2c\iota  t$ in $\cK$, such that 
$\mu \ast \eta t =\iota  \ast\pi $ and ${\widehat t}=\pi \ast\iota $. Then  
$({\widehat t}, {\widehat \mu }:=\pi \ast \mu  \ast\iota\iota , {\widehat \eta
  }:=\pi \ast \eta )$ is a monad in $\cK$. 
\end{itemize}
\end{lemma}

\begin{proof}
The proof is an easy computation using Definition \ref{def:pre_monad} of a
pre-monad and the properties $\iota $ and $\pi$ obey, so it is left to the
reader. 
\end{proof}

Improving \cite[Theorem 3.11]{WCP}, we obtain the following generalization of
a correspondence between monads in $\EM(\cK)$ and in $\cK$, observed by Lack
and Street in \cite{LS}.

\begin{theorem}\label{thm:composite}
For any monad $(k \1c t k, \mu ,\eta )$ and any 1-cell $k \1c s k$ in an
arbitrary 2-category $\cK$, there is a bijective correspondence between the
following structures. 
\begin{itemize}
\item[{(i)}] A monad $(t \1c {(s,\psi)} t, \nu , \vartheta)$ in $\EM^w(\cK)$;
\item[{(ii)}] A pre-monad $(st, \Theta, \vartheta)$ in $\cK$, such that 
\begin{equation}
\Theta \ast sts\mu  = s\mu  \ast \Theta t.
\label{eq:left_lin}
\end{equation}
\end{itemize}
\end{theorem}

\begin{proof}
The proof is built on the same constructions as \cite[Theorem 3.11]{WCP}.

Assume first that there exist 2-cells $\psi$ and $\nu$ as in part (i). A
multiplication $\Theta$ for the pre-monad in part (ii) is given by the same
formula of a `wreath product' in \cite[page 256]{LS}:
\begin{equation}\label{eq:M}
\Theta:= s\mu  \ast \nu t \ast ss\mu  \ast s\psi t.
\end{equation}
Its associativity is checked by the same computation as in the case of the
wreath product in \cite{LS}. By \eqref{eq:eta_lun} on one hand,
and by \eqref{eq:eta_2_a} and \eqref{eq:eta_run} on the other,
\begin{equation}\label{eq:M_eta}
\Theta \ast \vartheta st = s\mu  \ast \psi t \ast \eta st=\Theta \ast st
\vartheta, 
\end{equation}
proving \eqref{eq:preun_1}.
By applying \eqref{eq:eta_2_a} again, we conclude that $\Theta\ast \vartheta
\vartheta = \vartheta$, i.e. also \eqref{eq:preun_2} holds true. 
Condition \eqref{eq:norm} is proven by the following computation.
\begin{eqnarray*}
\Theta \ast \Theta st \ast \vartheta stst 
&=& s\mu  \ast \nu t \ast ss\mu  \ast s \psi t \ast s\mu st \ast \psi tst \ast
\eta stst \\ 
&=& s\mu  \ast \nu t \ast ss\mu  \ast ss\mu t \ast s \psi tt \ast st \psi t
\ast \psi tst \ast \eta stst \\ 
&=& s\mu  \ast s\mu t \ast s\mu  tt \ast \nu ttt \ast s \psi tt \ast \psi stt
\ast ts \psi t \ast \eta stst \\ 
&=& s\mu  \ast s\mu t \ast s\mu  tt \ast \psi ttt \ast t \nu tt \ast ts \psi t
\ast \eta stst \\
&=& s\mu  \ast s\mu t \ast \nu tt \ast s \psi t = \Theta.
\end{eqnarray*}
The second equality follows by \eqref{eq:psi_1}. 
The fourth and the fifth equalities 
follow by \eqref{eq:mu_2_a} and \eqref{eq:mu_2_b}, respectively.
Condition \eqref{eq:left_lin} follows by associativity of $\mu $. 
This proves that the data in part (i) determine a pre-monad as in part (ii).

Conversely, assume that there is a 2-cell $\Theta$ as in part (ii). The 2-cells
$\psi$ and $\nu$ in part (i) are constructed as
\begin{equation}\label{eq:psi_mu}
\psi:= \Theta \ast s\mu st \ast \vartheta tst \ast ts\eta  \qquad \qquad
\nu:= \Theta \ast s\eta s\eta .
\end{equation}
By \eqref{eq:left_lin}, 
\begin{equation}\label{eq:M_id}
s\mu  \ast \psi t = \Theta \ast s\mu st \ast \vartheta tst
\qquad \textrm{and}\qquad 
s\mu  \ast \nu t = \Theta \ast s\eta st.
\end{equation}
Moreover, by \eqref{eq:pre_associ}, \eqref{eq:left_lin} and
\eqref{eq:preun_1}, 
\begin{equation}\label{eq:M_psi}
\Theta\ast st \psi 
= \Theta \ast \Theta st \ast sts\mu st \ast st \vartheta tst \ast stts\eta 
= \Theta  \ast s\mu st \ast \Theta tst \ast \vartheta sttst \ast stts\eta .
\end{equation}
With identities \eqref{eq:M_id}, \eqref{eq:M_psi}, \eqref{eq:left_lin} and
\eqref{eq:preun_2} at hand, \eqref{eq:psi_1} is verified as
\begin{eqnarray*}
s\mu  \ast \psi t \ast t \psi 
&=& \Theta  \ast st \psi \ast s\mu ts \ast \vartheta tts
= \Theta  \ast s\mu st \ast \Theta tst \ast sts\mu tst \ast \vartheta \vartheta
ttst \ast tts\eta \\ 
&=& \Theta  \ast s\mu st \ast s\mu tst \ast \vartheta ttst \ast tts\eta 
= \psi \ast \mu s.
\end{eqnarray*}
Use next \eqref{eq:M_id}, \eqref{eq:M_psi}, \eqref{eq:left_lin} and
\eqref{eq:norm} to compute
\begin{equation}
s\mu  \ast \nu t \ast s\psi
=\Theta  \ast st \psi \ast s\eta ts
= \Theta  \ast \Theta st \ast \vartheta stst \ast sts\eta 
=\Theta  \ast sts\eta . \label{eq:m_mu_psi}
\end{equation}
In order to prove that \eqref{eq:mu_2_a} holds, apply \eqref{eq:m_mu_psi},
\eqref{eq:pre_associ} and \eqref{eq:M_id}:
\begin{equation}\label{eq:2.2}
s\mu  \ast \nu t \ast s \psi \ast \psi s 
= \Theta  \ast \psi st \ast tss\eta 
= \Theta  \ast \Theta st \ast s\mu stst \ast \vartheta tstst \ast ts\eta s\eta 
= s\mu  \ast \psi t \ast t\nu.
\end{equation}
Condition \eqref{eq:mu_2_b} is verified by comparing the last and third
expressions in \eqref{eq:2.2}, and using \eqref{eq:norm}:
$$
s\mu  \ast \psi t \ast \eta st \ast \nu
= \Theta  \ast \Theta  st \ast \vartheta stst \ast s \eta s\eta 
= \Theta  \ast s\eta s\eta 
=\nu.
$$
Condition \eqref{eq:eta_2_a} is proven by using \eqref{eq:M_id},
\eqref{eq:preun_1}, \eqref{eq:left_lin} and \eqref{eq:preun_2}:
$$
s\mu  \ast \psi t \ast t \vartheta 
= \Theta  \ast st\vartheta \ast s\mu  \ast \vartheta t 
= \Theta  \ast sts\mu  \ast \vartheta \vartheta t 
= s\mu  \ast \Theta t \ast \vartheta \vartheta t 
= s\mu  \ast \vartheta t.
$$
Condition \eqref{eq:mu_associ} follows by \eqref{eq:m_mu_psi},
\eqref{eq:pre_associ} and \eqref{eq:M_id}:
$$
s\mu  \ast \nu t \ast s\psi \ast \nu s
=  \Theta  \ast \Theta st \ast s\eta s\eta s\eta  
= \Theta  \ast s\eta st \ast s \Theta  \ast ss\eta s\eta 
= s\mu  \ast \nu t \ast s\nu.
$$
Condition \eqref{eq:eta_lun} is checked by applying \eqref{eq:m_mu_psi}:
\begin{equation}\label{eq:2.6}
s\mu  \ast \nu t \ast s\psi \ast \vartheta s
=  \Theta  \ast \vartheta st \ast s\eta 
=  \Theta  \ast s\mu st \ast \vartheta tst \ast \eta s\eta  
=\psi \ast \eta s.
\end{equation}
Finally, \eqref{eq:eta_run} is proven by making use of \eqref{eq:M_id},
\eqref{eq:preun_1} and comparing the second and last expressions in
\eqref{eq:2.6}:
$$
s\mu  \ast \nu t \ast s \vartheta
= \Theta  \ast st \vartheta \ast s\eta 
= \Theta  \ast \vartheta st \ast s\eta 
=\psi \ast \eta s.
$$
This proves that the data in part (ii) determine a monad as in part (i).

It remains to show that the above constructions are mutual inverses. Take 
2-cells $\nu$ and $\psi$ as in part (i). Use \eqref{eq:M} to associate a
2-cell $\Theta $ as in part (ii) to them, and then use \eqref{eq:psi_mu} to
define 2-cells $\nu'$ and $\psi'$ as in part (i) again. We obtain
\begin{eqnarray*}
\nu' 
&=& s\mu  \ast \nu t \ast s \psi  \ast s\eta s
= s\mu  \ast s\mu t \ast \nu tt \ast s \nu t \ast ss \vartheta
= s\mu  \ast \nu t \ast ss\mu  \ast s \psi t \ast st \vartheta \ast \nu\\
&=& s\mu  \ast s\mu t \ast \nu tt \ast s \vartheta t \ast \nu
= s\mu  \ast \psi t \ast \eta st \ast \nu 
=\nu.
\end{eqnarray*}
The first equality follows by unitality of $\mu $ and
the second equality follows by \eqref{eq:eta_run} and associativity of $\mu $.
The third equality is obtained by \eqref{eq:mu_associ} and
the fourth equality follows by \eqref{eq:eta_2_a} and associativity of $\mu $.
The penultimate equality is a consequence of \eqref{eq:eta_run} and the last
one follows by \eqref{eq:mu_2_b}. Also,
\begin{eqnarray*} 
\psi'
= s\mu  \ast \nu t \ast s \psi \ast s\mu s \ast \vartheta ts 
= s\mu  \ast s\mu t \ast \nu tt \ast s \psi t \ast \vartheta st \ast \psi
= s\mu  \ast \psi t \ast t \psi \ast \eta ts 
= \psi. 
\end{eqnarray*}
The first equality follows by unitality of $\mu $,
the second equality follows by \eqref{eq:psi_1} and associativity of $\mu $,
and the third equality is obtained by \eqref{eq:eta_lun}.
The last equality follows by \eqref{eq:psi_1} and by unitality of $\mu $.

In the opposite order, start with a 2-cell $\Theta $ as in part (ii). Apply
\eqref{eq:psi_mu} to construct 2-cells $\psi$ and $\nu$ as in part (i) and
then apply \eqref{eq:M} to obtain a 2-cell $\Theta '$ as in part (ii). It
satisfies  
\begin{eqnarray*} 
\Theta '
&=& s\mu  \ast \nu t \ast ss\mu  \ast s \psi t
= \Theta  \ast s\eta st \ast ss\mu \ast s \psi t
= s\mu  \ast \Theta t \ast st \psi t \ast s\eta tst\\
&=& \Theta  \ast s\mu st \ast \Theta tst \ast \vartheta sttst \ast s\eta tst
= \Theta  \ast \Theta st \ast \vartheta stst = \Theta .
\end{eqnarray*}
The second equality follows by the second identity in \eqref{eq:M_id}.
The third equality is obtained by \eqref{eq:left_lin}. 
The fourth equality is a consequence of \eqref{eq:M_psi} and unitality of $\mu
$. The penultimate equality follows by \eqref{eq:left_lin} and unitality of
$\mu $. The last equality is obtained by \eqref{eq:norm}.
\end{proof}

\begin{examples}
Examples of a composite pre-monad as in Theorem \ref{thm:composite}(ii) are
given, first 
of all, by wreath products in \cite{LS}. It is shown in \cite[Example 3.3]{LS}
that crossed products by Hopf algebras in \cite{Swe}, \cite{BCM}, \cite{DT}
are examples of a wreath product. As it was observed by Ross Street
\cite{priv} (see also \cite[Example 3.18]{WCP}), so are the crossed products
with coalgebras in \cite{Brz:crosp} and their generalizations to comonads in
\cite[Section 4.8]{Wis}.  

Examples of a composite pre-monad, which are not monads, are
provided by `weak smash products' in 
\cite[Section 3]{CaeDeG} (for a review see \cite[Example 3.16]{WCP}). This
includes smash products with weak bialgebras \cite{BNSz}.
Crossed products with weak bialgebras in \cite{RodRap} are also shown in
\cite[Section 4]{WCP} to provide examples.

Note, however, that (weak) crossed products with bialgebroids in \cite[Section
  4 \& Appendix]{BB} are not (pre-)monads of the kind in Theorem
\ref{thm:composite}(ii). 
Let $k$ be a commutative and associative unital ring.
A $k$-algebra $B$, measured by a left bialgebroid $H$ over a $k$-algebra $L$,
determines two monads, $(-)\otimes _k B$ on the category of $k$-modules and
$(-) \otimes_L B$ on the category of $L$-modules. 
In terms of the measuring $\cdot :H \otimes_k B \to B$ and the
comultiplication $h \mapsto \sum h_1\otimes_L h_2$ in $H$, consider the left
$L$-module map 
$$
H \otimes_k B \to B \otimes_L H,\qquad 
h \otimes_k b \mapsto \sum h_1\cdot b \otimes_L h_2.
$$
It equips the left $L$-module (or $L$-$k$ bimodule) $H$ with the structure of
a 1-cell $(-)\otimes_L B \1c {(-)\ot_L H} (-) \otimes_k B$ in
$\EM(\mathrm{CAT})$ (or in $\EM^w(\mathrm{CAT})$ if $\cdot$ is only a weak
measuring).  
However, if $L$ is a non-trivial $k$-algebra, this 1-cell has different source
and target.
Although in \cite{BB} the composite endofunctor $-\otimes_k B \otimes_L H$ on
the category of $k$-modules is proven to carry a monad structure, it is not a
composite of a monad with an endofunctor.
\end{examples}

\section{A pseudo-functor $\EM^w(\cK)\to \cK$} \label{sec:J^w}

Throughout this section (and the next one), we make two basic assumptions on
the 2-category $\cK$ we deal with:
\begin{itemize}
\item[{(i)}] Idempotent 2-cells in $\cK$ split;
\item[{(ii)}] $\cK$ admits Eilenberg-Moore constructions (EM constructions,
  for short) for monads. 
\end{itemize}
In more details, assumption (i) means that for any 2-cell $V \2c e V$ in
$\cK$, such that $e\ast e =e$, there exist a 1-cell ${\widehat V}$ and
2-cells ${\widehat V} \2c \iota  V$ and $V \2c \pi  {\widehat V}$, such that
$\pi\ast \iota = {\widehat V}$ and $\iota  \ast\pi  =e$. 
It is easy to see that the datum
$({\widehat V},\iota ,\pi )$ is unique up to an isomorphism ${\widehat V}
\2c \iota  V \2c {\pi '} {\widehat V}'$.  

Property (ii) of a 2-category was introduced by Street in \cite[page
151]{Street}. It means that the inclusion 2-functor $\cK \to \Mnd(\cK)$
(with underlying maps $k\mapsto (k \1c k k,k,k)$, $(k \1c V k') \mapsto (k\1c V
k', V \2c V V)$, $(V\2c \varrho W) \mapsto (V \2c \varrho W)$ on the
0-, 1-, and 2-cells, respectively) possesses a right 2-adjoint.
By \cite[Section 1]{LS}, property (ii) can be formulated equivalently by
saying that the inclusion 2-functor $I:\cK \to \EM(\cK)$ possesses a right
2-adjoint $J$.
Important properties of 2-categories admitting EM constructions for monads are
formulated in the following theorem. It is recalled from \cite[Theorem
  2]{Street} and \cite[Section 1]{LS}.

\begin{theorem}\label{thm:street}
In a 2-category $\cK$ which admits EM constructions for monads, any monad $(k
\1c t k,\mu ,\eta )$ determines an adjunction $(k \1c f J(t), J(t) \1c v k, k
\2c \eta  vf, fv \2c\epsilon  k)$ in $\cK$, such that $(t,\mu ,\eta
)=(vf,v\epsilon f,\eta )$. One can choose $JI=\cK$,
and the isomorphism corresponding to the 2-adjunction $(I,J)$ is given by the
mutually inverse functors 
\begin{eqnarray}\label{eq:2-adj}
&\cK(l,J(t)) \to \EM(\cK)(I(l),t), \ \ 
&(V \2c \omega W) \mapsto ((vV,v\epsilon V) \2c {v\omega} (vW,v\epsilon W))\\
&\EM(\cK)(I(l),t) \to \cK(l,J(t)), \ \ 
&((A,\alpha)\2c \varrho (B,\beta)) \mapsto (J(A,\alpha) \2c {J(\varrho)}
  J(B,\beta)), 
\nonumber
\end{eqnarray}
for any $0$-cell $l$ and monad $t$ in $\cK$.
\end{theorem}
The notations in Theorem \ref{thm:street} are used throughout, without further
explanation. 

\begin{lemma}\label{lem:psi_id}
Consider a 2-category $\cK$ which admits EM constructions for monads. For any
1-cell $(t,\mu ,\eta )\1c {(V,\psi)} (t',\mu ',\eta ')$ in $\EM^w(\cK)$, the
2-cell $Vv\epsilon \ast \psi v \ast \eta 'Vv:Vv \2c{} Vv$ in $\cK$ is
idempotent, and obeys the following identities. 
\begin{eqnarray}  
&& V\mu  \ast \psi t \ast \eta ' Vt \ast \psi = \psi;
\label{eq:l_proj}\\
&& Vv\epsilon  \ast \psi v \ast t' Vv\epsilon  \ast t' \psi v \ast t'\eta 'Vv=
Vv\epsilon  \ast \psi v. 
\label{eq:r_proj}
\end{eqnarray}
\end{lemma}

\begin{proof}
All statements follow easily by applying the interchange law,
\eqref{eq:1-cell} and unitality of the monad $t'$.
\end{proof}

The idempotent 2-cell in Lemma \ref{lem:psi_id}, corresponding to a 1-cell 
$(t,\mu ,\eta ) \1c {(V,\psi)} (t',\mu ',\eta ')$ 
in $\EM^w(\cK)$, is an identity 2-cell if and only if
    $\psi \ast \eta ' V = V\eta $, i.e.  $(V,\psi)$ is a 1-cell in $\EM(\cK)$.

Our next aim is to extend the 2-functor $J$ in Theorem \ref{thm:street} to
$\EM^w(\cK)$. Our method is reminiscent to the way $J$ is obtained from the
right adjoint of the inclusion 2-functor $\cK\to \Mnd(\cK)$.

\begin{lemma} \label{lem:1-cell}
Consider a 2-category $\cK$ which admits EM constructions for monads and in
which idempotent 2-cells split. 
For a 1-cell $(t,\mu ,\eta )\1c {(V,\psi)} (t',\mu ',\eta ')$ in $\EM^w(\cK)$,
denote a chosen splitting of the idempotent 2-cell in Lemma \ref{lem:psi_id}
by $Vv \2c \pi  {\widetilde V} \2c \iota Vv$. Then $({\widetilde V},
{\widetilde \psi}:=\pi\ast Vv\epsilon  \ast \psi v \ast t' \iota)$ is a 1-cell
$IJ(t) \1c {} t'$ in $\EM(\cK)$. 
\end{lemma}

\begin{proof}
By the interchange law, ${\widetilde \psi}\ast \eta ' {\widetilde V}=\pi \ast
\iota \ast\pi  \ast \iota = {\widetilde V}$. Furthermore, by \eqref{eq:r_proj}
and 
\eqref{eq:1-cell},
\begin{eqnarray*}
{\widetilde \psi} \ast t' {\widetilde \psi}
&=&\pi  \ast Vv\epsilon  \ast \psi v \ast t' Vv\epsilon  \ast t' \psi v \ast
t' t' \iota   
=\pi  \ast Vv\epsilon  \ast V\mu v \ast \psi tv \ast t'\psi v \ast t' t' \iota
\\ 
&=&\pi  \ast Vv\epsilon  \ast \psi v \ast \mu ' Vv\ast t' t' \iota 
= {\widetilde \psi} \ast \mu ' {\widetilde V}.  
\end{eqnarray*}
\end{proof}

\begin{lemma} \label{lem:2-cell}
Consider a 2-category $\cK$ which admits EM constructions for monads and in
which idempotent 2-cells split. 
For any 2-cell $(V,\psi) \2c \varrho (W,\phi)$ in $\EM^w(\cK)$,
${\widetilde \varrho}:=\pi\ast Wv\epsilon  \ast \varrho v \ast \iota$ is a
2-cell in $\EM(\cK)$, between the 1-cells $({\widetilde V},{\widetilde \psi})$
and $({\widetilde W},{\widetilde \phi})$ in Lemma \ref{lem:1-cell} (where
$\pi$ and $\iota $ denote chosen splittings of both idempotent
2-cells in Lemma \ref{lem:psi_id}, corresponding to the 1-cells $(V,\psi)$ and
$(W,\phi)$).   
\end{lemma}

\begin{proof}
Apply \eqref{eq:l_proj} (in the first equality), \eqref{eq:2-cell_a} (in the
third equality) and \eqref{eq:2-cell_b} (in the penultimate equality) to
conclude that  
\begin{eqnarray*}
{\widetilde \varrho} \ast {\widetilde \psi}
&=&\pi  \ast Wv\epsilon  \ast \varrho v \ast Vv\epsilon  \ast \psi v \ast t'
\iota   
=\pi  \ast Wv\epsilon  \ast W\mu v \ast \varrho tv \ast \psi v \ast t' \iota \\
&=&\pi  \ast Wv\epsilon  \ast W\mu v \ast \phi tv \ast t'\varrho v\ast t'
\iota  
=\pi  \ast Wv\epsilon  \ast \phi v \ast t' \iota  \ast t'\pi  \ast
t'Wv\epsilon \ast t' \varrho v \ast
t' \iota  \\
&=& {\widetilde \phi}\ast t' {\widetilde \varrho}.
\end{eqnarray*}
\end{proof}

\begin{theorem}\label{thm:J^w}
Consider a 2-category $\cK$ which admits EM constructions for monads and in
which idempotent 2-cells split. 
The following maps determine a pseudo-functor $J^w:\EM^w(\cK)\to \cK$.
\begin{itemize}
\item[]
For a \underline{0-cell} $t$, $J^w(t):= J(t)$.
\item[]
For a \underline{1-cell} $t \1c {(V,\psi)} t'$,
  $J^w(V,\psi):=J({\widetilde V},{\widetilde \psi})$, where the 1-cell
  $({\widetilde V},{\widetilde \psi})$ in $\EM(\cK)$ is described in Lemma
  \ref{lem:1-cell}. 
  That is, denoting by $\iota $ and $\pi$ a chosen splitting of
  the idempotent 2-cell in Lemma \ref{lem:psi_id}, 
  $J^w(V,\psi)$ is the unique 1-cell $J^w(t) \1c {} J^w(t')$ in $\cK$ for which
  $v'\epsilon ' J^w(V,\psi) =\pi \ast Vv\epsilon  \ast \psi v \ast t'\iota $.
\item[] For a \underline{2-cell} $(V,\psi) \2c \varrho (W,\phi)$,
  $J^w(\varrho):= J({\widetilde \varrho})$, where the 2-cell ${\widetilde
  \varrho}$ in $\EM(\cK)$ is described in Lemma \ref{lem:2-cell}. That is,
  $J^w({\varrho})$ is the unique 2-cell $J^w(V,\psi)\2c {}
  J^w(W,\phi)$ in $\cK$ for which 
  $v' J^w({\varrho})=\pi \ast Wv\epsilon  \ast \varrho v \ast \iota$.
\end{itemize}
The pseudo-natural isomorphism class of $J^w$ is independent of
  the choice of the 2-cells $\iota$ and $\pi$ in its construction.
\end{theorem}

\begin{proof}
Let us fix splittings $(\pi,\iota)$ of the idempotent 2-cells in Lemma
\ref{lem:psi_id}, for all 1-cells $(V,\psi)$ in $\EM^w(\cK)$. 

By construction, 
$v' J^w(\psi \ast \eta ' V) =\pi  \ast \iota \ast\pi  \ast \iota = {\widetilde
  V}$, for any 1-cell $t\1c {(V,\psi)} t'$ in $\EM^w(\cK)$. Hence
$J^w$ preserves unit 2-cells $(V,\psi)\2c {\psi \ast \eta ' V} (V,\psi)$. 
For 2-cells $(V,\psi) \2c \varrho (W,\phi) \2c \tau (U,\theta)$ in
$\EM^w(\cK)$, it follows by \eqref{eq:2-cell_b} (applied to $\varrho$) that 
\begin{eqnarray*}
v' J^w(\tau)\ast v' J^w(\varrho)
&=&\pi  \ast Uv\epsilon  \ast \tau v \ast Wv\epsilon  \ast \varrho v \ast
\iota\\ 
&=&\pi  \ast Uv\epsilon  \ast U \mu v \ast \tau tv \ast \varrho v \ast \iota 
= v' J^w(\tau \bullet \varrho).
\end{eqnarray*}
We conclude by the isomorphism \eqref{eq:2-adj} that $J^w$ preserves the
vertical composition. 

For an identity 1-cell $t \1c{(k,t)} t$, the idempotent 2-cell in Lemma
\ref{lem:psi_id} is the identity 2-cell $v$ by the adjunction relation
$v\epsilon  \ast \eta v =v$. Hence any splitting of it yields
mutually inverse isomorphisms $v\2c{\pi_k} {\widetilde k}$ and ${\widetilde k}
\2c {\iota_k} v$. They give rise to an isomorphism $J^w(t)=J(v,v\epsilon) \2c 
{J(\pi_k)} J^w(k,t)=J({\widetilde k}, \pi_k \ast v \epsilon\ast t \iota_k)$ with
the inverse $J(\iota_k)$. Thus $J^w$ preserves identity 1-cells up to
isomorphism. (In particular, we can choose for the definition of $J^w$ a
trivial splitting $v \2c v v \2c v v$, in which case the 1-cell $({\widetilde
k}, {\widetilde t})$ in Lemma \ref{lem:1-cell} is equal to $(v,v\epsilon
)$. Applying the isomorphism \eqref{eq:2-adj}, we conclude that with this
choice, $J^w$ strictly preserves identity 1-cells,
i.e. $J^w(k,t)=J(v,v\epsilon)=J^w(t)$.) 

In order to investigate the preservation of the horizontal composition,
consider different splittings $(\pi,\iota )$ and $(\pi',\iota ')$ of the
idempotent 2-cell in Lemma \ref{lem:psi_id}, for some
1-cell $(V,\psi)$, and denote the corresponding 1-cells in Lemma
\ref{lem:1-cell} by $(\widetilde{V},\widetilde{\psi})$ and
$(\widetilde{V}',\widetilde{\psi}')$, respectively. Applying \eqref{eq:r_proj}
(for $(\pi',\iota ')$) and \eqref{eq:l_proj} (for $(\pi,\iota )$),
$$
\pi' \ast Vv\epsilon  \ast \psi v \ast t'\iota ' \ast t'\pi ' \ast t'\iota 
=\pi ' \ast Vv\epsilon  \ast \psi v \ast t' \iota 
=\pi' \ast \iota \ast\pi  \ast Vv\epsilon  \ast \psi v \ast t'\iota .
$$
Hence $(\widetilde{V},\widetilde{\psi})
\2c{\pi '\ast \iota} (\widetilde{V}',\widetilde{\psi}')$
is an iso 2-cell in $\EM(\cK)$, so $J(\widetilde{V},\widetilde{\psi})
\2c{J(\pi'\ast \iota)} J(\widetilde{V}',\widetilde{\psi}')$
is an iso 2-cell in $\cK$. 

For 1-cells $(V,\psi),(W,\phi):t \1c {} t'$ and
$(V',\psi'),(W',\phi'):t' \1c {} t''$ and 2-cells $(V,\psi) \2c \varrho
(W,\phi)$ and $(V',\psi') \2c \varrho' (W',\phi')$ in $\EM^w(\cK)$, the
idempotent 2-cell in Lemma \ref{lem:psi_id} corresponding to the 1-cell
$(V',\psi') \circ (V,\psi)$ comes out as $V'Vv\epsilon  \ast V' \psi v \ast
\psi' Vv \ast \eta '' V' V v$. We claim that 
it has a splitting given by the mono 2-cell $V'\iota  \ast \iota'
J^w(V,\psi)$ and the epi 
2-cell $\pi' J^w(V,\psi) \ast V'\pi $, where $(\pi,\iota )$ and $(\pi',\iota
')$ are the chosen splittings 
of the idempotent 2-cells in Lemma \ref{lem:psi_id}, corresponding to the
1-cells $(V,\psi)$ and $(V',\psi')$ in $\EM^w(\cK)$, respectively. Indeed, by
construction of $J^w$ (its action on a 1-cell $(V,\psi)$), 
\eqref{eq:l_proj} and \eqref{eq:r_proj},   
\begin{eqnarray} \label{eq:splittings}
V'\iota  \ast \iota' J^w(V,\psi) &\ast&\pi ' J^w(V,\psi) \ast V'\pi \nonumber
\\ 
&=& V'\iota  \nast V'\pi \nast V' Vv\epsilon  \nast V' \psi v \nast V' t' \iota
\nast V' t'\pi  \nast \psi' Vv \nast \eta '' V'Vv \nonumber \\
&=& V'Vv\epsilon  \ast V'\psi v \ast \psi' Vv \ast \eta '' V'Vv. 
\end{eqnarray}
Denote by $V'Vv \2c {\pi _2} {\widetilde{V'V}} \2c {\iota_2} V'Vv$ 
the canonical splitting of this idempotent which was chosen to define
$J^w$ on the 1-cell $(V',\psi')\circ (V,\psi)=(V'V,V'\psi \ast \psi' V)$. By
considerations in the previous paragraph, there are mutually inverse iso
2-cells  
$J(\pi ' J^w(V,\psi) \ast V'\pi \ast \iota_2):
J^w(V'V,V'\psi \ast \psi' V) \2c {} J^w(V',\psi')J^w(V,\psi)$
and
$J(\pi _2 \ast V'\iota  \ast \iota 'J^w(V,\psi)): J^w(V',\psi')J^w(V,\psi)\2c
{} J^w(V'V,V'\psi \ast \psi' V)$. 
In order to see their naturality, observe that
$$
v''J^w(\varrho'\circ \varrho)
= \pi _2 \ast W'Wv\epsilon  \ast W'W\mu v \ast W' \varrho tv \ast W'
\psi v \ast \varrho' Vv\ast \iota_2.
$$
On the other hand,
\begin{eqnarray*}
&v''& J^w(\varrho')J^w(\varrho)\\
&=&\pi ' J^w(W,\phi) \nast W'\pi  \nast W' Wv\epsilon  \nast W' \phi v \nast
  \varrho' Wv \nast  V' \iota  \nast V'\pi  \nast V' W v\epsilon  \nast V'
  \varrho v \nast V' \iota  \nast \iota' J^w (V,\psi)\\ 
&=&\pi ' J^w(W,\phi)\nast W'\pi  \nast W' W v\epsilon  \nast W' W \mu v \nast W'
  \phi tv \nast W' t' \varrho v \ast \varrho' Vv \nast V' \iota  \ast \iota'
  J^w(V,\psi) \\ 
&=&\pi ' J^w(W,\phi) \nast W'\pi  \nast W'Wv\epsilon  \nast W'W\mu v \nast W'
  \varrho tv \nast W'\psi v \ast \varrho' Vv\nast V'\iota  \ast \iota'
  J^w(V,\psi). 
\end{eqnarray*}
The second equality follows by applying \eqref{eq:2-cell_b}, and the third one
follows by applying \eqref{eq:2-cell_a}, for $\varrho$.
With this information in mind, we conclude that 
\begin{eqnarray*}
v'' J^w(\varrho' \circ \varrho)\ast
&\pi_2& \ast V'\iota  \ast \iota' J^w(V,\psi) \\
&=&\pi _2 \ast W'Wv\epsilon  \ast W'W\mu v \ast W' \varrho tv \ast W' \psi v
\ast \varrho' Vv \ast V'\iota  \ast \iota' J^w(V,\psi)\\
&=&\pi _2 \ast W'\iota  \ast \iota' J^w(W,\phi)\ast v'' J^w(\varrho')
J^w(\varrho).  
\end{eqnarray*}
Thus naturality of $J(\pi_2 \ast V'\iota \ast \iota' J^w(V,\psi))$
follows by the isomorphism \eqref{eq:2-adj}.
It remains to check its associativity and unitality.
For a further 1-cell $(V'',\psi''):t'' \to t'''$ in $\EM^w(\cK)$, use the
notation  
$V''V'Vv \2c {\pi _3} 
{\widetilde{V''V'V}}
\2c {\iota_3} V''V'Vv$ for the canonically split idempotent in the
construction of $J^w$ on $(V'',\psi'')\circ (V',\psi')\circ
(V,\psi)=(V''V'V,V''V'\psi \ast V'' \psi' V \ast \psi'' V'V)$. By
\eqref{eq:splittings}, the associativity condition  
\begin{eqnarray*}
&&\pi _3 \ast V'' \iota_2 \ast \iota'' J^w(V'V,V'\psi \ast \psi'V) \ast v'''
  J^w(V'',\psi'')J(\pi_2\ast V'\iota \ast \iota'J^w(V,\psi)) \\ 
&=&\pi _3 \ast V'' \iota_2 \ast V''\pi_2 \ast V''V'\iota \ast V''\iota
  'J^w(V,\psi) \ast \iota'' J^w(V',\psi')J^w(V,\psi) \\
&=&\pi _3 \ast V'' V' \iota  \ast V'' \iota ' J^w(V,\psi) \ast \iota''
  J^w(V',\psi')J^w(V,\psi)\\
&=&\pi _3 \ast V'' V' \iota  \ast \iota_2' J^w(V,\psi) \ast\pi _2'
  J^w(V,\psi)\ast V'' \iota ' J^w(V,\psi) \ast \iota'' J^w(V',\psi')J^w(V,\psi) 
\end{eqnarray*}
holds. The 2-cells $\iota_k$ and $\pi_k$, splitting the idempotent
(identity) 2-cell in Lemma \ref{lem:psi_id} corresponding to a unit 1-cell
$(k,t)$, are mutual inverses. Hence also the unitality conditions
\begin{eqnarray*}
&&v'J(\pi \ast V \iota_k \ast \iota J^w(k,t))\ast v' J^w(V,\psi) J(\pi_k) =
\pi \ast V \iota_k \ast V \pi_k \ast \iota = 
v' J^w(V,\psi)\\
&&v' J(\iota_{k'} J^w(V,\psi)) \ast v' J(\pi_{k'}) J^w(V,\psi) =
\iota_{k'} J^w(V,\psi) \ast \pi_{k'} J^w(V,\psi) 
 = v' J^w(V,\psi)
\end{eqnarray*}
hold. Thus we conclude by the isomorphism \eqref{eq:2-adj} that $J^w$ 
preserves also the horizontal composition up to a coherent
family of iso 2-cells, i.e. that it is a pseudo-functor.

Finally, we investigate the ambiguity of the pseudo-functor $J^w$, caused by a
free choice of the splittings of the idempotent 2-cells in Lemma
\ref{lem:psi_id}. Take two collections $\{(\pi,\iota)\}$ and
$\{(\pi',\iota')\}$ of splittings (indexed by the 1-cells in
$\EM^w(\cK)$). The pseudo-functors $J^w$ and $J'^w$, associated to
both families of splittings, are pseudo-naturally isomorphic via
$J^w(t)=J'^w(t)$ and $J^w(V,\psi) \2c {J(\pi'\ast \iota)} J'^w(V,\psi)$, for
any 0-cell $t$ and 1-cell $(V,\psi)$ in $\EM^w(\cK)$. 
\end{proof}

Consider a 2-category $\cK$ which admits EM constructions for monads and in
which idempotent 2-cells split. 
We can regard a 1-cell $t \1c {(V,\psi)} t'$ in $\EM(\cK)$ as a 1-cell in 
$\EM^w(\cK)$. Choosing a trivial splitting $Vv \2c{Vv} Vv\2c {Vv} Vv$ of the
identity 2-cell, the corresponding 1-cell $IJ(t) \1c {({\widetilde{V}},
{\widetilde{\psi}})} t'$ in Lemma \ref{lem:1-cell} comes out as the 1-cell
$(Vv, Vv\epsilon  \ast \psi v)$ in $\EM(\cK)$. By 1-naturality of the counit of
the 2-adjunction $(I,J)$, we have
$(v' J(V,\psi), v'\epsilon ' J(V,\psi))= (Vv, Vv\epsilon  \ast \psi v)$.
From this and from the isomorphism \eqref{eq:2-adj} it follows that
$$
J^w(V,\psi)=J(Vv, Vv\epsilon  \ast \psi v)=J(V,\psi).
$$
Similarly, we can regard a 2-cell $(V,\psi)\2c \varrho (W,\phi)$ in $\EM(\cK)$
as a 2-cell in $\EM^w(\cK)$. The corresponding 2-cell $({\widetilde{V}},
{\widetilde{\psi}}) \2c {\widetilde \varrho} ({\widetilde{W}},
{\widetilde{\phi}})$ in Lemma \ref{lem:2-cell} is equal to the 2-cell
$Wv\epsilon  \ast \varrho v:(Vv, Vv\epsilon  \ast \psi v) \2c {} (Wv,
Wv\epsilon  \ast \phi v)$ in $\EM(\cK)$. By the 2-naturality condition 
$v' J(\varrho)=Wv\epsilon  \ast \varrho v$ 
and the isomorphism \eqref{eq:2-adj} we obtain that
$$
J^w(\varrho)=J(Wv\epsilon  \ast \varrho v)=J(\varrho).
$$
Summarizing, we proved that the pseudo-functor $J^w:\EM^w(\cK)\to \cK$ in
Theorem \ref{thm:J^w} can be chosen such that the 2-functor $J:\EM(\cK)\to
\cK$ in Theorem \ref{thm:street} factorizes through the obvious inclusion
$\EM(\cK)\hookrightarrow \EM^w(\cK)$ and $J^w$. 

The pseudo-functor $J^w$ in Theorem \ref{thm:J^w} takes a monad
$((s,\psi),\nu,\vartheta)$ in $\EM^w(\cK)$ to a monad $J^w(s,\psi)$ in $\cK$,
with multiplication $J^w(s,\psi)J^w(s,\psi) \2c \cong J^w((s,\psi) \circ
(s,\psi)) \2c {J^w(\nu)} J^w(s,\psi)$ and unit $J^w(t) \2c \cong J^w(k,t) \2c
{J^w(\vartheta)} J^w(s,\psi)$. Applying to this monad in $\cK$ a hom 2-functor
$\cK(l,-):\cK \to \mathrm{CAT}$ (for any 0-cell $l$ in $\cK$), we obtain a
monad in $\mathrm{CAT}$.  
Our next aim is to describe its Eilenberg-Moore category.  

\begin{lemma}\label{lem:AC_obj}
Consider a 2-category $\cK$ which admits EM constructions for monads and in
which idempotent 2-cells split. 
Let $l$ be a 0-cell and $(k \1c t k,\mu ,\eta )$ be a monad in $\cK$ and 
let $t \1c {(s,\psi)} t$ be a 1-cell in $\EM^w(\cK)$.
There is a bijective correspondence between the following
  structures. 
\begin{itemize}
\item[{(i)}] Pairs $(l \1c V J^w(t), J^w(s,\psi)V \2c \zeta V)$, consisting of
  a 1-cell $V$ and a 2-cell $\zeta$ in $\cK$; 
\item[{(ii)}] Pairs $((l \1c W k,tW \2c \varrho W), sW \2c \lambda W)$,
  consisting of a 1-cell $I(l) \1c {(W,\varrho)} t$ in $\EM(\cK)$ and 
(regarding $(W,\varrho)$ as a 1-cell in $\EM^w(\cK)$), a 2-cell
$(s,\psi)\circ (W,\varrho) \2c \lambda (W,\varrho)$ in $\EM^w(\cK)$.
\end{itemize}
\end{lemma}

\begin{proof}
Denote by $\iota $ and $\pi$ the splitting of the idempotent 2-cell in Lemma
\ref{lem:psi_id}, corresponding to the 1-cell $(s,\psi)$ in $\EM^w(\cK)$,
that was chosen to construct $J^w(s, \psi)$. For the 1-cell
$(W,\varrho)$ in $\EM^w(\cK)$, choose the trivial splitting of the identity
  2-cell $W \2c W W$, so that $J^w(W,\varrho)= J(W,\varrho)$.

By \eqref{eq:2-adj}, there is a bijection between the 1-cells 
$I(l) \1c {(W,\varrho)} t$ in $\EM(\cK)$ as in part (ii), and the 1-cells
$V:=J(W,\varrho)=J^w(W,\varrho):l \1c {} J(t)=J^w(t)$ in $\cK$ as in part
(i). 
In order to see that it extends to the stated bijection, 
take first a 2-cell 
$\lambda$ in $\EM^w(\cK)$ as in part (ii). Then there 
is a 2-cell $\zeta:= \big(J^w(s,\psi)V \2c \cong J^w((s,\psi)\circ
(W,\varrho)) \2c{J^w(\lambda)}  V\big)$ in $\cK$ as in part (i).  
Conversely, for a 2-cell 
$\zeta$ in $\cK$ as in part (i), put
$\lambda:=v \zeta \ast\pi V$. It satisfies
\begin{eqnarray} 
\varrho \ast t \lambda
&=& v\epsilon V \ast tv \zeta \ast t\pi V
= v \zeta \ast v\epsilon  J^w(s,\psi)V \ast t\pi V\nonumber\\
&=& v \zeta \ast\pi V \ast sv\epsilon  V \ast \psi v V \ast t\iota V \ast t\pi
V 
= \lambda \ast s \varrho \ast \psi W. 
\label{eq:psi_mod}\label{eq:c_psi}
\end{eqnarray}
The second equality follows by the interchange law and the third one follows
by construction of the pseudo-functor $J^w$, cf. Theorem \ref{thm:J^w}.
The last equality follows by \eqref{eq:r_proj}. Together with the unitality of
$\varrho$, this proves that $\lambda$ is a 2-cell in $\EM^w(\cK)$, as needed.

The above two constructions can be seen to be mutual inverses. Take first a
pair $(V,\zeta)$ as in part (i) and iterate both constructions. The result
is 
$\big(V, J^w(s,\psi)V \2c \cong J^w((s,\psi)\circ (W,\varrho))\2c {J^w(v \zeta
  \ast \pi v)} V\big) =(V,\zeta)$.
In the opposite order, a datum $((W,\varrho),\lambda)$ is taken to
$\big((W,\varrho), sW \2c{\pi V} v J^w(s,\psi)V \2c \cong v J^w((s,\psi)\circ
(W,\varrho)) \2c {v J(\lambda)} W\big)=((W,\varrho),\lambda \ast \iota V
\ast\pi V)$. 
The resulting 2-cell $\lambda \ast \iota V \ast\pi V$ in $\cK$ is equal to
$\lambda$ since by \eqref{eq:psi_mod} and unitality of $\varrho$,
\begin{equation} \label{eq:la_ip}
\lambda \ast \iota V \ast\pi V 
= \lambda \ast s \varrho \ast \psi W \ast \eta sW
= \varrho \ast t \lambda \ast \eta sW
=\lambda.
\end{equation}
\end{proof}

The following extends \cite[Proposition 3.1]{LS} and also \cite[Theorem
  3.4]{CaeDeG}. 

\begin{proposition}\label{prop:modules}
Consider a 2-category $\cK$ which admits EM constructions for monads and in
which idempotent 2-cells split. 
Let $l$ be a 0-cell and $(k \1c t k,\mu ,\eta )$ be a monad in $\cK$ and 
let $(t \1c {(s,\psi)} t, \nu, \vartheta)$ be a monad in $\EM^w(\cK)$.
The following categories are isomorphic.
\begin{itemize}
\item[{(i)}] The Eilenberg-Moore category
$\A$ of the monad 
$\cK(l,J^w(s,\psi))
:\cK(l,J^w(t))$ $ \to \cK(l,J^w(t))$; 
\item[{(ii)}] The Eilenberg-Moore category ${\B}$ of the monad
$\cK(l,{\widehat{st}}): \cK(l,k) \to \cK(l,k)$, where the monad
  ${\widehat{st}}$ is obtained from the pre-monad $st$ in Theorem
  \ref{thm:composite} in the way described in Lemma \ref{lem:pre_monad};
\item[{(iii)}] The category ${\mathcal C}$, with
\item[] \underline{objects} that are pairs $((W,\varrho), \lambda)$,
  consisting of a 1-cell $I(l) \1c {(W,\varrho)} t$ in $\EM(\cK)$ and a 2-cell  
$(s,\psi)\circ (W,\varrho)\2c \lambda (W,\varrho)$ in  $\EM^w(\cK)$, 
satisfying
\begin{eqnarray}
&& \lambda \bullet ((s,\psi)\circ \lambda) = \lambda \bullet (\nu \circ
  (W,\varrho));\label{eq:lam_associ}\\
&& (W,\varrho)=\lambda \bullet (\vartheta \circ (W,\varrho));\label{eq:lam_uni}
\end{eqnarray}
\item[]\underline{morphisms} $((W,\varrho), \lambda)\to
  ((W',\varrho'),\lambda')$ that are 2-cells $(W,\varrho) \2c \alpha
  (W',\varrho')$ in $\EM(\cK)$ such that  
\begin{equation}\label{eq:C_mor}
\lambda' \bullet ((s,\psi) \circ \alpha)= \alpha \bullet \lambda.
\end{equation}
\end{itemize}
\end{proposition}

\begin{proof}
Denote by $\iota $ and $\pi$ the splitting of the idempotent 2-cell in Lemma
\ref{lem:psi_id}, corresponding to the 1-cell $(s,\psi)$ in $\EM^w(\cK)$,
that was chosen to construct $J^w(s,\psi)$.
For the 1-cell
$(W,\varrho)$ in $\EM^w(\cK)$, choose the trivial splitting of the identity
  2-cell $W \2c W W$, so that $J^w(W,\varrho)= J(W,\varrho)$.
Introduce shorthand notations ${\overline s}:=J^w(s,\psi)$ and
$V:=J^w(W,\varrho)$.  

Isomorphism of \underline{${\A}$ and $ {\mathcal C}$}.
In light of Lemma \ref{lem:AC_obj}, any object in ${\A}$ is of the
form $(J^w(W,\varrho), {\overline s}V \2c \cong J^w((s,\psi)\circ (W,\varrho))
\2c {J^w(\lambda)} V)$, for a unique 1-cell $I(l)\1c{(W,\varrho)} t$ in
$\EM(\cK)$ and a unique 2-cell $(s,\psi)\circ (W,\varrho)
\2c{\lambda}(W,\varrho)$ in $\EM^w(\cK)$. So we only need to show that
$\lambda$ satisfies \eqref{eq:lam_associ}, i.e. the equality 
\begin{equation}\label{eq:c_mu}
\lambda \ast s\lambda =\lambda \ast s\varrho \ast \nu W
\end{equation}
of 2-cells in $\cK$, if and only if ${\overline s}V \2c \cong
J^w((s,\psi)\circ (W,\varrho)) \2c{J^w(\lambda)}V$ is an
associative action, and $\lambda$ satisfies \eqref{eq:lam_uni}, i.e. 
\begin{equation}\label{eq:c_eta}
W=\lambda \ast s \varrho\ast \vartheta W
\end{equation}
if and only if this ${\overline s}$-action on $V$ is unital.
Compose the associativity condition 
$
J^w(\lambda)\ast  J^w((s,\psi)\circ \lambda) = J^w(\lambda) \ast J^w(\nu \circ
(W,\varrho))
$ 
with $v$ horizontally on the left, and compose it with 
the chosen split epimorphism $ssW \to vJ^w(ssW, ss\varrho \ast
  s\psi W \ast \psi sW)$
on the right (i.e. on the `top').
It yields the equality
$$
\lambda \ast s \lambda \ast \big(ss \varrho \ast s \psi W \ast \psi s
W \ast \eta ss W\big)
= \lambda \ast s \varrho \ast \nu W \ast \big(ss \varrho \ast s \psi W \ast
\psi s W \ast \eta ss W\big).
$$
Making use of \eqref{eq:psi_mod}, the left hand side is easily shown to be
equal to $\lambda \ast s \lambda$. As far as the right hand side is
concerned, use associativity of $\varrho$ (in the first equality),
\eqref{eq:mu_2_a} and \eqref{eq:mu_2_b} (in the second and third equalities,
respectively) to see that 
\begin{eqnarray*}
\lambda \nast s \varrho \nast \nu W\nast ss \varrho \nast  s \psi W \nast \psi
s W \nast \eta ssW 
&=& \lambda \nast s \varrho \nast s\mu  W \nast \nu t W \nast s \psi W \nast
\psi s W \nast \eta ss W\\  
= \lambda \ast s\varrho \ast s\mu  W \ast \psi t W \ast t \nu W \ast \eta ss W
&=& \lambda \ast s\varrho \ast \nu W.
\end{eqnarray*}
This proves that $J^w(\lambda)$ is associative if and only if
\eqref{eq:c_mu} holds.
Similarly, the unitality condition $J^w(\lambda) \ast J^w(\vartheta\circ
(W,\varrho))  = V$ is 
equivalent to $\lambda \ast \iota V \ast\pi V \ast s \varrho\ast \vartheta W
=W$, hence by \eqref{eq:la_ip} it is equivalent to \eqref{eq:c_eta}.  
Thus the bijection in Lemma \ref{lem:AC_obj} restricts to a bijection
between the objects in ${\A}$ and the objects in ${\mathcal C}$.

For a 2-cell $(W,\varrho) \2c \alpha (W',\varrho')$ in $\EM(\cK)$, the
condition $J^w(\lambda')\ast J^w((s,\psi)\circ \alpha) = J^w(\alpha) \ast
J^w(\lambda)$ 
(expressing that $J^w(\alpha)=J(\alpha)$ is a morphism in ${\A}$)
is equivalent to 
$\lambda' \ast s\alpha \ast \iota V \ast \pi V=
\alpha\ast \lambda \ast \iota V \ast \pi V$.
The right hand side is equal to $\alpha\ast \lambda$ by \eqref{eq:la_ip} and
the left hand side is equal to 
$$
\lambda' \ast s\alpha \ast s\varrho \ast \psi W \ast \eta sW
= \lambda' \ast s\varrho' \ast \psi W' \ast \eta sW' \ast s\alpha
= \varrho' \ast t \lambda' \ast \eta sW' \ast s\alpha
= \lambda' \ast s\alpha,
$$
using that $\alpha$ is a 2-cell in $\EM(\cK)$, \eqref{eq:psi_mod} and
unitality of $\varrho'$. Hence $J^w(\alpha)=J(\alpha)$ is a morphism in ${\A}$
if and only if \eqref{eq:C_mor} holds. 
Thus we conclude by the isomorphism \eqref{eq:2-adj} that the 2-functor $J$
induces an (obviously functorial) isomorphism between the
morphisms in ${\A}$ and the morphisms in ${\mathcal C}$. 

Isomorphism of \underline{${\B}$ and $ {\mathcal C}$}.
In view of \eqref{eq:M_eta}, we can choose
${\widehat{st}}= v {\overline s } f$ as 1-cells in $\cK$. Moreover, taking
axioms \eqref{eq:pre_associ}, \eqref{eq:preun_1} and \eqref{eq:norm} of a
pre-monad into account, \begin{equation}\label{eq:M_norm}
\iota f \ast\pi f \ast \Theta  
=\Theta  \ast \iota fst \ast\pi fst
=\Theta  \ast st\iota f \ast st\pi f 
=\Theta .
\end{equation}
For an object $(l \1c W k, v {\overline s} fW \2c \gamma W)$ in ${\B}$, put
\begin{equation}\label{eq:rho_lambda}
\varrho:= \gamma \ast\pi fW \ast s\mu W \ast \vartheta tW
\qquad \textrm{and}\qquad 
\lambda:=\gamma \ast\pi fW \ast s\eta W.
\end{equation}
We show that $((W,\varrho),\lambda)$ is an object in ${\mathcal C}$. Recall
that associativity and unitality of $\gamma$ read as  
$$
\gamma \ast v {\overline s} f \gamma 
= \gamma \ast\pi fW \ast \Theta  W \ast \iota f\iota fW 
\qquad \textrm{and}\qquad 
W=\gamma \ast\pi fW \ast \vartheta W,
$$
respectively. Hence using associativity of $\gamma$ and \eqref{eq:M_norm} (in
the second equality) and applying the first identity in \eqref{eq:M_id} (in the
last equality), 
\begin{eqnarray}
\varrho \ast t \gamma \ast t\pi fW
&=& \gamma \ast v {\overline s}f \gamma \ast\pi f\pi fW \ast s\mu st W \ast
\vartheta tst W 
\nonumber\\
&=& \gamma \ast\pi fW \ast \Theta W \ast s\mu stW \ast \vartheta tstW
= \gamma \ast\pi fW \ast s\mu W \ast \psi tW.
\label{eq:rho_id}
\end{eqnarray}
Moreover, apply associativity of $\gamma$ and \eqref{eq:M_norm} (in the
second equality) and use \eqref{eq:left_lin} (in the third equality)
to obtain 
\begin{eqnarray}
\lambda \ast s\varrho
&=& \gamma \ast v {\overline s}f \gamma \ast\pi f\pi fW \ast s\eta stW \ast
ss\mu  W \ast s \vartheta tW \nonumber \\
&=&\gamma \ast\pi fW \ast \Theta W \ast sts\mu W \ast st \vartheta tW \ast
s\eta tW \nonumber \\ 
&=& \gamma \ast\pi fW \ast s\mu W \ast \Theta tW \ast st \vartheta t W \ast
s\eta t W  
= \gamma \ast\pi fW.
\label{eq:lambda_id}
\end{eqnarray}
In the last equality we used \eqref{eq:M_eta} and
that (since $\mu =v\epsilon f$) the interchange law yields 
$s\mu  \ast \iota ft \ast\pi ft = \iota f  \ast\pi f \ast s\mu $.
With these identities at hand, associativity of $\varrho$ is checked as
\begin{eqnarray*}
\varrho \ast t \varrho
&=& \varrho \ast t \gamma \ast t\pi fW \ast ts \mu  W \ast t\vartheta t W
= \gamma \ast\pi fW \ast s\mu W \ast s\mu tW \ast \psi tt W \ast t \vartheta t
W\\ 
&=& \gamma \ast\pi fW \ast s\mu W \ast s\mu tW \ast \vartheta ttW 
= \varrho \ast \mu W.
\end{eqnarray*}
The second equality follows by \eqref{eq:rho_id} and by associativity of
$\mu $. 
In the third equality we applied \eqref{eq:eta_2_a}.
The last equality follows by associativity of $\mu $ and the form of $\varrho$
in \eqref{eq:rho_lambda}.
The unitality condition  $\varrho\ast \eta W =W$ follows by unitality of $\mu $
and unitality of $\gamma$. Conditions \eqref{eq:c_psi}, \eqref{eq:c_mu}
and \eqref{eq:c_eta} are proven by
\begin{eqnarray*}
\varrho \ast t \lambda 
&=& \varrho \ast t\gamma \ast t\pi fW \ast ts\eta W
= \gamma \ast\pi fW \ast s\mu W \ast \psi t W \ast ts\eta W
= \lambda \ast s\varrho \ast \psi W;\\
\lambda \ast s\lambda 
&=& \gamma \ast v {\overline s} f \gamma \ast\pi f\pi fW \ast s\eta s\eta  W 
=\gamma \ast\pi fW \ast \nu W 
= \lambda \ast s \varrho \ast \nu W;\\
W&=& \gamma \ast\pi fW \ast \vartheta W = \lambda \ast s \varrho \ast
\vartheta W. 
\end{eqnarray*}
In each case, the last equality follows by \eqref{eq:lambda_id}. 
In the first computation, the second equality follows by \eqref{eq:rho_id}.
In the second equality of the second computation we used 
associativity of $\gamma$ together with \eqref{eq:M_norm} and we applied the
expression of $\nu$ in \eqref{eq:psi_mu}.
In the first equality of the last computation we used unitality of
$\gamma$. 
This proves that $((W,\varrho),\lambda)$ is an object in ${\mathcal C}$.

Conversely, for an object $((W,\varrho),\lambda)$ in ${\mathcal C}$, put
\begin{equation}\label{eq:gamma}
\gamma:= \lambda \ast s \varrho \ast \iota fW.
\end{equation}
It is associative as
\begin{eqnarray*}
\gamma \ast\pi fW \ast \Theta W \ast \iota f\iota fW
&=& \lambda \ast s \varrho \ast \Theta W \ast \iota f\iota fW\\
&=& \lambda \ast s \varrho \ast \nu W \ast ss \varrho \ast s \psi W \ast sts
\varrho \ast \iota f\iota fW\\
&=& \lambda \ast s \lambda \ast ss\varrho \ast s \psi W \ast sts \varrho \ast
\iota f\iota fW\\ 
&=& \lambda \ast s \varrho \ast st \lambda \ast sts \varrho \ast \iota f\iota
fW 
= \gamma \ast v {\overline s} f \gamma. 
\end{eqnarray*}
The first equality follows by \eqref{eq:gamma} and \eqref{eq:M_norm}.
In the second equality we substituted $\Theta $ by its expression in
\eqref{eq:M} and we used associativity of $\varrho$ twice.
In the third equality we applied \eqref{eq:c_mu} and in the fourth one we used
\eqref{eq:c_psi}. 
By \eqref{eq:eta_2_a} and unitality of $\mu $, $\iota f \ast\pi f \ast
\vartheta = s\mu  \ast \psi t \ast \eta st \ast \vartheta = \vartheta$. Hence
the unitality condition $\gamma \ast\pi fW \ast \vartheta W =W$ follows by
\eqref{eq:c_eta}. 
This proves that $(W,\gamma)$ is an object in ${\B}$.

Let us see that the above constructions are mutual inverses.
Starting with an object $(W,\gamma)$ of ${\B}$ and iterating the above
constructions, we re-obtain $(W,\gamma)$ by \eqref{eq:lambda_id}.
In the opposite order, applying both constructions to an object
$((W,\varrho),\lambda)$ of ${\mathcal C}$, we obtain
$((W, \lambda \ast s\varrho \ast \iota fW \ast\pi fW \ast s\mu W \ast \vartheta
tW), \lambda \ast s \varrho \ast \iota fW \ast\pi fW \ast s\eta W)$. 
Since $\iota f \ast\pi f \ast s\mu  = s\mu  \ast \iota ft \ast\pi ft = s\mu
\ast \Theta  t\ast \vartheta stt$, 
axiom \eqref{eq:preun_2} of a pre-monad, associativity of $\varrho$ and
\eqref{eq:c_eta} imply that   
$$
\lambda \ast s\varrho \ast \iota fW \ast\pi fW \ast s\mu W \ast \vartheta tW
= \lambda \ast s \varrho \ast s\mu W \ast \vartheta tW 
= \lambda \ast s \varrho \ast \vartheta W \ast \varrho
= \varrho.
$$
Also, by \eqref{eq:M_eta},
unitality of $\mu $, \eqref{eq:c_psi} and unitality of $\varrho$,
$$
\lambda \ast s \varrho \ast \iota fW \ast\pi fW \ast s\eta W
= \lambda \ast s \varrho \ast s\mu W \ast \psi tW \ast \eta s\eta W
= \lambda \ast s \varrho \ast \psi W \ast \eta s W
= \lambda.
$$
Hence we constructed a bijection between the objects of ${\B}$ and
${\mathcal C}$. 
It is immediate by the form of the bijection between the objects that a 2-cell
$W \2c \alpha W'$ in $\cK$ is a morphism in ${\B}$ if and only if it
is a morphism in ${\mathcal C}$.
\end{proof}

\section{Weak liftings} \label{sec:lift}

If $\cK$ is a 2-category which admits EM constructions for monads, 
then `liftings' of 1-, and 2-cells for monads in $\cK$ arise as images under
the right 2-adjoint $J$ of the inclusion 2-functor $\cK \to \EM(\cK)$, see
\cite{PowWat}. In this section we discuss `weak' liftings and the role what
the pseudo-functor $J^w$ plays in their description.

\begin{definition}\label{def:1-cell_lift}
Consider a 2-category $\cK$ which admits EM constructions for monads.
We say that a 1-cell $k \1c V k'$ in $\cK$ possesses a {\em weak lifting} for
some monads $(k \1c t k, \mu ,\eta )$ and $(k' \1c {t'} k', \mu ', \eta ')$ in
$\cK$ if there exist a 1-cell $J(t) \1c {\overline V} J(t')$ and a split mono
2-cell $v' {\overline V} \2c \iota Vv$ (with retraction denoted by $Vv \2c \pi
v'{\overline V}$).  
\end{definition}

If in a 2-category $\cK$ which admits EM constructions for monads, also
idempotent 2-cells split, then we know by Theorem \ref{thm:J^w} that,
for every 1-cell $t \1c {(V,\psi)} t'$ in $\EM^w(\cK)$, the underlying 
1-cell $k\1c V k'$ in $\cK$ possesses a weak lifting $J^w(V,\psi)$ for $t$ and
$t'$. As we will see later in this section, in fact, in such a
2-category $\cK$, up-to an isomorphism, every weak lifting arises in this
way. This extends assertions about 1-cells in \cite[Lemma 3.9 and Theorem
3.10]{PowWat}.  

\begin{definition}\label{def:2-cell_lift}
Consider a 2-category $\cK$ which admits EM constructions for monads.
Let $(k \1c t k,\mu ,\eta )$ and $(k' \1c {t'} k',\mu ',\eta ')$ be monads,
and $k \1c V k'$ and $k \1c W k'$ be 1-cells in $\cK$, such that there exist
their weak liftings $(J(t) \1c {\overline V} J(t'),\iota_V,\pi_V)$ and 
$(J(t) \1c {\overline W} J(t'),\iota_W,\pi_W)$ for $t$ and $t'$. 
For a 2-cell $V \2c \omega W$ in $\cK$, we say that 
\begin{itemize}
\item[] a 2-cell 
${\overline V} \2c {\ilift \omega} {\overline W}$ 
is a {\em weak $\iota $-lifting} of $\omega$ if 
$\iota _W \ast v' {\ilift \omega}= \omega v \ast \iota_V$;
\item[] a 2-cell 
${\overline V} \2c {\plift \omega} {\overline W}$ 
is a   {\em weak $\pi$-lifting} of $\omega$ if 
$v' {\plift \omega} \ast\pi _V =\pi _W \ast \omega v$.
\end{itemize}
\end{definition}

Throughout, indices of $\iota $ and $\pi$ are omitted, as they can be
reconstructed without ambiguity from the context.

By the isomorphism \eqref{eq:2-adj}, the weak $\iota$-lifting or weak
$\pi$-lifting of a 2-cell is unique, provided that it exists.
Moreover, if a 2-cell $\omega$ in Definition \ref{def:2-cell_lift} possesses
both a weak $\iota$-lifting $\ilift \omega$ and a weak $\pi$-lifting $\plift
\omega$, then 
$$
v' \ilift \omega 
=\pi \ast \iota \ast v' \ilift \omega 
=\pi  \ast \omega v \ast \iota 
= v' \plift \omega \ast\pi  \ast \iota
= v' \plift \omega.
$$
Hence in view of the isomorphism \eqref{eq:2-adj}, $\ilift \omega = \plift
\omega$. 

Proposition \ref{prop:2-cell_lift} below, about weak liftings of 2-cells in a
(nice enough) 2-category, extends statements about 2-cells in \cite[Lemma 3.9
and Theorem 10]{PowWat}. Therein, notions and notations introduced in
Corollary \ref{cor:Mnd} are used.

\begin{proposition}\label{prop:2-cell_lift}
Consider a 2-category $\cK$ which admits EM constructions for monads and in
which idempotent 2-cells split. 
Let $t\1c {(V,\psi)} t'$ and $t \1c {(W,\phi)} t'$ be 1-cells in $\EM^w(\cK)$
and $V \2c \omega W$ be a 2-cell in $\cK$. Denote by $\iota $ and $\pi$
the splittings of both idempotent 2-cells in Lemma \ref{lem:psi_id},
corresponding to the 1-cells $(V,\psi)$ and $(W,\phi)$, that were chosen to
construct $J^w(V,\psi)$ and $J^w(W,\phi)$, respectively. 

\begin{itemize}
\item[{(1)}] The following assertions are equivalent.
\begin{itemize}
\item[{(i)}] $\omega$ is a 2-cell $(V,\psi)\2c{} (W,\phi)$ in $\Mnd^\iota
  (\cK)$; 
\item[{(ii)}] $\omega t \ast \psi \ast \eta ' V:(V,\psi) \2c {} (W,\phi)$ is a
  2-cell in $\EM^w(\cK)$;
\item[{(iii)}] $\pi \ast \omega v \ast \iota:(v' J^w(V,\psi), v'\epsilon
  'J^w(V,\psi)) \2c 
  {} (v' J^w(W,\phi), v'\epsilon 'J^w(W,\phi))$ is a 2-cell in $\EM(\cK)$ such
  that 
  $\iota \ast\pi  \ast \omega v \ast \iota = \omega v \ast \iota$;
\item[{(iv)}] $\omega$ possesses a weak $\iota $-lifting ${\ilift \omega}:
  (J^w(V,\psi), \iota,\pi) \to (J^w(W,\phi), \iota,\pi)$.
\end{itemize}
If these equivalent statements hold, then 
$v' J^wG^\iota (\omega)
=\pi \ast \omega v \ast \iota$, i.e. 
$\ilift \omega = 
J^wG^\iota (\omega)
$.

\item[{(2)}] The following assertions are equivalent.
\begin{itemize}
\item[{(i)}] $\omega$ is a 2-cell $(V,\psi)\2c{} (W,\phi)$ in $\Mnd^\pi (\cK)$;
\item[{(ii)}] $\phi\ast \eta ' W \ast \omega:(V,\psi) \2c {} (W,\phi)$ is a
  2-cell in $\EM^w(\cK)$;
\item[{(iii)}] $\pi \ast \omega v \ast \iota:(v' J^w(V,\psi), v'\epsilon
  'J^w(V,\psi)) \2c {} (v' J^w(W,\phi), v'\epsilon 'J^w(W,\phi))$ is a 2-cell
  in $\EM(\cK)$ such that $\pi \ast \omega v \ast \iota \ast\pi =\pi \ast
  \omega v$; 
\item[{(iv)}] $\omega$ possesses a weak $\pi$-lifting $\plift
  \omega:(J^w(V,\psi), \iota,\pi) \to (J^w(W,\phi), \iota,\pi)$. 
\end{itemize}
If these equivalent statements hold, then 
$v' J^w G^\pi (\omega)
=\pi \ast \omega v \ast \iota$, i.e. 
$\plift \omega = 
J^wG^\pi (\omega)
$.

\item[{(3)}] The following assertions are equivalent.
\begin{itemize}
\item[{(i)}] $\phi \ast t' \omega = \omega t \ast \psi$; 
\item[{(ii)}] $\phi\ast \eta ' W \ast \omega$ and $\omega t \ast \psi \ast
  \eta ' V$ 
  are (necessarily equal) 2-cells $(V,\psi) \2c {} (W,\phi)$ in $\EM^w(\cK)$; 
\item[{(iii)}] $\pi \ast \omega v \ast \iota:(v' J^w(V,\psi), v'\epsilon
  'J^w(V,\psi)) \2c 
  {} (v' J^w(W,\phi), v'\epsilon 'J^w(W,\phi))$ is a 2-cell in $\EM(\cK)$ such
  that 
  $\iota \ast\pi  \ast \omega v = \omega v \ast \iota \ast\pi $; 
\item[{(iv)}] $\omega$ possesses both a weak $\iota$-lifting and a weak
  $\pi$-lifting 2-cell $(J^w(V,\psi), \iota,\pi) \to (J^w(W,\phi), \iota,\pi)$
  (which are necessarily equal).
\end{itemize}
\end{itemize}
\end{proposition}

\begin{proof}
\underline{(1)} \underline{(i)$\Leftrightarrow$(ii)} 
This equivalence follows by Lemma \ref{lem:2-cells}(1).

\underline{(ii)$\Leftrightarrow$(iii)} 
The 2-cell $\pi \ast \omega v \ast \iota$ in $\cK$ is a 2-cell in $\EM(\cK)$
if and only if $v'\epsilon ' J^w(W,\phi) \ast t'\pi  \ast t' \omega v \ast t'
\iota  =\pi  \ast \omega v \ast \iota \ast v'\epsilon ' J^w(V,\psi)$.
Compose this equality by $f$ horizontally on the right, and compose it
vertically by $\iota f$ on the left and by $t'\pi f\ast t'V\eta $ on the right. 
By virtue of \eqref{eq:l_proj}, \eqref{eq:r_proj} and the adjunction relation
$\epsilon f\ast f\eta =f$, the resulting equivalent condition is equivalent to
\eqref{eq:omega_a}.  
Property $\iota \ast\pi  \ast \omega v \ast \iota = \omega v \ast \iota$ is
equivalent to $\iota \ast\pi  \ast \omega v \ast \iota \ast\pi  = \omega v
\ast \iota\ast\pi $, what is easily 
seen to be equivalent to \eqref{eq:omega_b}.
Thus we conclude by Lemma \ref{lem:2-cells}(1) (i)$\Leftrightarrow$(iii). 

\underline{(iii)$\Leftrightarrow$(iv)} By the isomorphism \eqref{eq:2-adj},
$\pi\ast \omega v \ast \iota$ is a 2-cell in $\EM(\cK)$ if and only if there
is a 2-cell $J^w(V,\psi)\2c {\ilift \omega} J^w(W,\phi)$ 
in $\cK$ such that $v' {\ilift \omega} =\pi \ast \omega v \ast \iota$.
Clearly, $\iota \ast\pi  \ast \omega v \ast \iota = \omega v \ast \iota$ if
and only if $\ilift \omega$ is a weak $\iota$-lifting of $\omega$.

If the equivalent statements (i)-(iv) hold, then 
$$
v'J^w(\omega t \ast \psi \ast \eta ' V) 
=\pi \ast Wv\epsilon  \ast \omega t v \ast \psi v \ast \eta ' Vv \ast \iota
=\pi \ast \omega v \ast \iota \ast\pi  \ast \iota 
=\pi  \ast \omega v \ast \iota.
$$

\underline{(2)} \underline{(i)$\Leftrightarrow$(ii)} 
This equivalence follows by Lemma \ref{lem:2-cells}(2).

\underline{(ii)$\Leftrightarrow$(iii)} 
As we have seen in the proof of part (1), $\pi \ast \omega v \ast \iota$ is a
2-cell in $\EM(\cK)$ if and only if \eqref{eq:omega_a} holds. Property  
$\pi \ast \omega v \ast \iota \ast\pi  =\pi  \ast \omega v$ is equivalent to 
$\iota  \ast\pi  \ast \omega v \ast \iota \ast\pi  = \iota \ast\pi  \ast
\omega v$ hence to the first condition in Lemma
\ref{lem:2-cells}(2)(iii). Thus we conclude by Lemma
\ref{lem:2-cells}(2) (i)$\Leftrightarrow$(iii).  

\underline{(iii)$\Leftrightarrow$(iv)} is proven by the same reasoning as in
part (1).

If the equivalent statements (i)-(iv) hold, then  
$$
v' J^w(\phi \ast \eta ' W \ast \omega) 
=\pi \ast Wv\epsilon  \ast \phi v \ast \eta ' Wv \ast \omega v \ast \iota =\pi
\ast \iota \ast\pi  \ast \omega v \ast \iota =\pi  \ast \omega v \ast \iota. 
$$

\underline{(3)} These equivalences follow immediately by Lemma
\ref{lem:2-cells}(3) and parts (1) and (2) in the current theorem. 
\end{proof}

For suggesting the following theorem, the author is grateful
  to the referee.

Consider a 2-category $\cK$ which admits EM constructions for monads. To any
monads $(k \1c t k,\mu ,\eta )$ and $(k' \1c {t'} k',\mu ',\eta ')$ in $\cK$,
we can associate categories $\mathrm{Lift}^\iota (t,t')$ and
$\mathrm{Lift}^\pi (t,t')$, as follows. In both categories objects are
quadruples $(V, \overline{V},\iota ,\pi )$ such that the 1-cell $J(t) \1c
{\overline V} J(t')$ in $\cK$ is a weak lifting of the 1-cell $k \1c V k'$,
corresponding to the split monic 2-cell $v' \overline{V} \2c \iota  Vv$, with
retraction $Vv \2c \pi v'\overline{V}$. Morphisms $(V, \overline{V},\iota ,\pi
) \to (W, \overline{W},\iota ,\pi )$ are pairs of 2-cells $(\omega, {\overline
  \omega})$ in $\cK$ such that $\overline{V} \2c {\overline \omega}
\overline{W}$ is a weak $\iota$-lifting, respectively, a weak $\pi$-lifting,
of $V \2c \omega W$. Composition of morphisms is defined via component-wise
composition of 2-cells in $\cK$.  

\begin{theorem} \label{thm:lifteq}
For any 2-category $\cK$ which admits EM constructions for monads and in
which idempotent 2-cells split, and for any monads $(k \1c t k,\mu ,\eta )$
and $(k' \1c {t'} k',\mu ',\eta ')$ in $\cK$, the following assertions hold.
\begin{itemize}
\item[{(1)}] $\mathrm{Lift}^\iota (t,t')$ is equivalent to the category
  $\Mnd^\iota(\cK) (t,t')$. 
\item[{(2)}] $\mathrm{Lift}^\pi (t,t')$ is equivalent to the category
  $\Mnd^\pi(\cK) (t,t')$. 
\end{itemize}
For $t=t'$, these equivalences are also strong monoidal, with respect to the
monoidal structure of $\mathrm{Lift}^{\iota/\pi} (t,t)$ induced by the
horizontal composition in $\cK$.
\end{theorem}

\begin{proof}
{\underline {(1)}}
For any 1-cell $t \1c {(V,\psi)} t'$ in $\EM^w(\cK)$, denote by 
$Vv \2c {\pi _c} v' J^w(V,\psi) \2c {\iota_c} Vv $ the chosen splitting of
the idempotent 2-cell in Lemma \ref{lem:psi_id}, used to construct
$J^w(V,\psi)$. 

By Corollary \ref{cor:Mnd} and Theorem \ref{thm:J^w}, there is a
pseudo-functor $J^wG^\iota :\Mnd^\iota(\cK)\to \cK$. By Proposition
\ref{prop:2-cell_lift}(1) (i)$\Rightarrow$(iv), it induces a functor  
$G: \Mnd^\iota(\cK) (t,t') \to \mathrm{Lift}^\iota (t,t')$,
with object map $(V,\psi) \mapsto (V, J^w(V,\psi), \iota_c,\pi_c)$ and
morphism map $\omega \mapsto (\omega, J^wG^\iota (\omega))$.   

In the opposite direction, consider the functor $F:\mathrm{Lift}^\iota (t,t')
\to \Mnd^\iota(\cK) (t,t')$, with the object map  
\begin{equation}\label{eq:lift_psi}
(V, \overline{V},\iota ,\pi ) \mapsto (V, \psi:= \iota f \ast v'\epsilon '
  {\overline V} f \ast t'  \pi f \ast t' V\eta : t'V \2c {} Vt) 
\end{equation}
and morphism map $(\omega,{\overline \omega})\mapsto \omega$. 
By the form of $\psi$ in \eqref{eq:lift_psi} and the adjunction relation $v
\epsilon\ast \eta v=v$, it follows that  
\begin{equation}\label{eq:psi_ip}
Vv\epsilon \ast \psi v = \iota \ast v' \epsilon' {\overline V} \ast t' \pi,
\end{equation}
thus in particular $\iota_c \ast \pi_c = Vv\epsilon \ast \psi v \ast \eta' Vv=
\iota \ast \pi$.
Using \eqref{eq:psi_ip} together with the form of $\psi$ in
\eqref{eq:lift_psi} and naturality, it is easily checked that $\psi$ satisfies
\eqref{eq:1-cell}, i.e. $(V,\psi)$ is an object in $\Mnd^\iota(\cK)(t,t')$.
Applying \eqref{eq:psi_ip} together with \eqref{eq:l_proj} and
\eqref{eq:r_proj}, respectively, we conclude that
$$
v'\epsilon 'J^w(V,\psi)\ast t'\pi_c \ast t'\iota  =\pi _c \ast \iota \ast
v'\epsilon ' {\overline  V}\quad \textrm{and}\quad 
\pi\ast \iota_c \ast v'\epsilon 'J^w(V,\psi)= v'\epsilon ' {\overline V}\ast
t'\pi \ast t' \iota_c. 
$$
That is, there are 2-cells $(v' J^w(V,\psi),
v'\epsilon' J^w(V,\psi)) \2c {\pi \ast \iota_c  } (v' {\overline V},
v'\epsilon'{\overline V})$ and $(v' {\overline V}, v'\epsilon'{\overline V})
\2c {\pi_c \ast \iota} (v' J^w(V,\psi), v'\epsilon' J^w(V,\psi))$ in
$\EM(\cK)$. By \eqref{eq:2-adj} they induce mutually inverse isomorphisms
$J^w(V,\psi) \2c {J( \pi \ast \iota_c)} {\overline V}$ and ${\overline V} \2c
{J(\pi_c \ast \iota)} J^w(V,\psi)$ in $\cK$. Both of these 2-cells are weak
$\iota$-liftings of the identity 2-cell $V \2c V V$. Hence, for any morphism
$(V, {\overline V},\iota,\pi) \1c {(\omega, {\overline \omega})} (W,
{\overline W},\iota,\pi)$ in $\mathrm{Lift}^\iota (t,t')$, the composite 
$J(\pi_c \ast \iota) \ast {\overline \omega} \ast J(\pi \ast \iota_c 
):J^w(V,\psi) \2c {} J^w(W,\phi)$ is a weak $\iota$-lifting of $\omega$ (where
both $\psi$ and $\phi$ are defined via \eqref{eq:lift_psi}).
Thus it follows by Proposition \ref{prop:2-cell_lift}(1) (iv)$\Rightarrow$(i)
that $\omega$ is a morphism $(V,\psi)\to (W,\phi)$ in
$\Mnd^\iota(\cK)(t,t')$. This proves that $F$ is a well defined functor.  

For any object $(V,\psi)$ in $\Mnd^\iota(\cK) (t,t')$, we obtain
$FG(V,\psi)=(V,\psi)$ by \eqref{eq:l_proj}, \eqref{eq:r_proj} and unitality of
$\mu $. Evidently, also $FG(\omega)=\omega$. For any object $(V,{\overline
  V},\iota ,\pi )$ of $\mathrm{Lift}^\iota (t,t')$, we obtain $GF(V,{\overline
  V},\iota ,\pi ) = (V,J^w(V,\psi),\iota_c,\pi_c)$. 
The mutually inverse isomorphisms 
$(V,{\overline V}, \iota, \pi) \1c {(V,J(\pi_c \ast \iota))}
 (V,J^w(V,\psi),\iota_c,\pi_c)$ and $(V,J^w(V,\psi),\iota_c,\pi_c) \1c
{(V,J(\pi \ast \iota_c ))} (V,{\overline V}, \iota, \pi)$
in $\mathrm{Lift}^\iota (t,t')$ define, in turn, mutually inverse natural
isomorphisms between the identity functor and $GF$. Indeed, for any morphism
$(V,{\overline V},\iota ,\pi )  \1c{(\omega,{\overline \omega})} (W,{\overline
  W},\iota ,\pi )$ in $\mathrm{Lift}^\iota (t,t')$, we conclude by Corollary
and \ref{cor:Mnd}, Theorem \ref{thm:J^w} that 
$$
v'J^w G^\iota(\omega) \ast\pi _c \ast \iota 
=\pi _c \ast Wv\epsilon  \ast \omega tv \ast \psi v \ast \eta 'Vv \ast \iota 
=\pi _c \ast \omega v \ast \iota 
=\pi _c \ast \iota \ast v' {\overline \omega}.
$$
Hence naturality follows by the isomorphism \eqref{eq:2-adj}.

It remains to prove strong monoidality of $G$ in the $t=t'$ case. 
Recall from the proof of Theorem \ref{thm:J^w} that the coherent natural
isomorphisms
$J^w(V',\psi') J^w(V,\psi) \2c {j_{V',V}} J^w((V',\psi') \circ (V,\psi))$
and $J^w(t) \2c {j_0} J^w(k,t)$,
rendering $J^w$ (hence $J^w G^\iota$) a pseudo-functor, arise as
weak $\iota$-liftings of identity 2-cells, for any 1-cells $(V',\psi'),
(V,\psi): t \to t$ in $\EM^w(\cK)$. Hence they induce a strong monoidal
structure $G(V',\psi') G(V,\psi) \2c {(V'V, j_{V',V})} G((V',\psi') \circ
(V,\psi))$ and $(k,J^w(t)) \2c {(k,j_0)} G(k,t)$ of $G$.

Part {\underline{(2)}} is proven symmetrically.
\end{proof}

\section{Applications}\label{sec:appli}

In this section we collect from the literature several situations where weak
liftings occur.
The following corollary is a consequence of Proposition \ref{prop:2-cell_lift}
and Theorem \ref{thm:lifteq}. 

\begin{corollary}\label{cor:comonad_lift}
Consider a 2-category $\cK$ which admits EM constructions for monads and in
which idempotent 2-cells split. 
Let $(k\1c t k, \mu ,\eta )$ be a monad and $(k\1c c k,\delta ,\varepsilon )$
be a comonad in $\cK$ . 
\begin{itemize}
\item[{(1)}] The following assertions are equivalent.
\begin{itemize}
\item[{(i)}] There exists a comonad $((c,\psi),\delta ,\varepsilon )$ in
  $\Mnd^\iota (\cK)$. That is, there exists a 1-cell $t \1c {(c,\psi)} t$ in
  $\EM^w(\cK)$, satisfying 
\begin{eqnarray}
&& \delta t \ast \psi 
= cc\mu  \ast c \psi t \ast \psi ct \ast t\delta t \ast t \psi \ast t\eta c;
\label{eq:d_i}\\
&& \varepsilon t \ast \psi 
= \mu  \ast t\varepsilon t \ast t \psi \ast t\eta c.
\label{eq:e_i}
\end{eqnarray}
 \item[{(ii)}] There is a comonad
$(t \1c {(c,\psi)} t, \delta t \ast \psi \ast \eta c, \varepsilon t \ast \psi
   \ast \eta c)$ in $\EM^w(\cK)$.
\item[{(iii)}] There is a comonad $(J^w(t) \1c {\overline c} J^w(t), \ilift
  \delta , \ilift \varepsilon )$ in $\cK$ such that $\ilift \delta $ is a weak
  $\iota$-lifting of $\delta $ and $\ilift \varepsilon $ is a weak
  $\iota$-lifting of $\varepsilon $. 
\end{itemize}
If these equivalent statements hold, then we say shortly that the comonad
$(\overline c, \ilift \delta , \ilift \varepsilon )$ is a weak $\iota$-lifting
of the comonad $(c,\delta ,\varepsilon )$ for the monad $(t,\mu ,\eta )$.

\item[{(2)}] The following assertions are equivalent.
\begin{itemize}
\item[{(i)}] There exists a comonad $((c,\psi),\delta ,\varepsilon )$ in
  $\Mnd^\pi (\cK)$. That is, there exists a 1-cell $t \1c {(c,\psi)} t$ in
  $\EM^w(\cK)$, satisfying 
\begin{eqnarray}
&& c\psi \ast \psi c \ast t\delta   
= cc\mu  \ast c \psi t \ast \psi ct \ast \eta cct \ast \delta t \ast \psi;
\label{eq:d_p}\\
&& t\varepsilon =\varepsilon t \ast \psi. 
\label{eq:e_p}
\end{eqnarray}
 \item[{(ii)}] There is a comonad
$(t \1c {(c,\psi)} t, c\psi \ast \psi c \ast \eta cc \ast \delta , \eta  \ast
   \varepsilon )$ in $\EM^w(\cK)$.
\item[{(iii)}] There is a comonad $(J^w(t) \1c {\overline c} J^w(t), \plift
  \delta , \plift \varepsilon )$ in $\cK$ such that $\plift \delta $ is a weak
  $\pi$-lifting of $\delta $ and $\plift \varepsilon $ is a weak $\pi$-lifting
  of $\varepsilon $. 
\end{itemize}
If these equivalent statements hold, then we say shortly that the comonad
$(\overline c, \plift \delta , \plift \varepsilon )$ is a weak $\pi$-lifting
of the comonad $(c,\delta ,\varepsilon )$ for the monad $(t,\mu ,\eta )$.

\item[{(3)}] The following assertions are equivalent.
\begin{itemize}
\item[{(i)}] There exists a 1-cell $t \1c {(c,\psi)} t$ in $\EM^w(\cK)$,
  satisfying 
\begin{eqnarray}
&& c\psi \ast \psi c \ast t\delta   
= \delta t \ast \psi;
\label{eq:d_b}\\
&& t\varepsilon =\varepsilon t \ast \psi. 
\label{eq:e_b}\label{eq:SK_14}
\end{eqnarray}
\item[{(ii)}] There is a comonad $(J^w(t) \1c {\overline c} J^w(t), \overline
  \delta , \overline \varepsilon )$ in $\cK$ such that $\overline \delta $ is
  both a weak $\iota$-lifting and a weak $\pi$-lifting of $\delta $ and
  $\overline \varepsilon $ is both a weak $\iota$-lifting and a weak
  $\pi$-lifting of $\varepsilon $.  
\end{itemize}   
\end{itemize}
\end{corollary}

Note that by Lemma \ref{lem:2-cells} and Lemma \ref{lem:lift_comp}, in parts (1)
and (2) of Corollary \ref{cor:comonad_lift}, assertions (i) and (ii) are
equivalent in case of an arbitrary 2-category $\cK$.

Let us stress the (tiny) difference between a 2-cell $tc \2c \psi ct$ in $\cK$
occurring in Corollary \ref{cor:comonad_lift}(3)(i), and a mixed distributive
law. A 2-cell $\psi$ in Corollary \ref{cor:comonad_lift}(3)(i) satisfies three
of the identities defining a mixed distributive law: compatibility with the
multiplication of the monad (as $(c,\psi)$ is a 1-cell in $\EM^w(\cK)$),
compatibility with the comultiplication of the comonad (by \eqref{eq:d_b}) and
compatibility with the counit of the comonad (by \eqref{eq:e_b}). However, the
fourth condition on a mixed distributive law, compatibility $\psi \ast \eta c
= c\eta $ with the unit of the monad, does not appear in Corollary
\ref{cor:comonad_lift}(3)(i) -- it plays no role in a weak lifting.

\begin{example}\label{ex:weak_entw}
Generalizing a mixed distributive law of a monad and a comonad (in particular
in the bicategory $\mathrm{BIM}$),
weak entwining structures were introduced by Caenepeel and De
Groot in \cite{CaeDeG}. The axioms are obtained by weakening the compatibility
conditions of a mixed distributive law with the unit of the monad and the
counit of the comonad.
Precisely, a {\em weak entwining structure} in an
arbitrary 2-category $\cK$ consists of a monad $(k\1c t k,\mu ,\eta )$, a
comonad $(k\1c c k,\delta ,\varepsilon )$ and a 2-cell $tc \2c \psi ct$
subject to the following conditions.
\begin{eqnarray}
&& \psi \ast \mu c 
= c\mu  \ast \psi t \ast t\psi;
\label{eq:wentw_m} \label{eq:SK_9}\\
&& \delta t \ast \psi 
= c\psi \ast \psi c \ast t\delta ;
\label{eq:wentw_d}\\
&& \psi \ast \eta c 
= c\varepsilon t \ast c\psi \ast c\eta c \ast \delta ;
\label{eq:wentw_u}\\
&& \varepsilon t \ast \psi 
= \mu  \ast t\varepsilon t \ast t \psi \ast t\eta c.
\label{eq:wentw_e}
\end{eqnarray}
We claim that under these assumptions $((c,\psi),\delta t \ast \psi \ast \eta
c, \varepsilon t \ast \psi \ast \eta c)$ is a comonad in $\EM^w(\cK)$.
For that, we need to show that axioms \eqref{eq:wentw_m}-\eqref{eq:wentw_e}
imply \eqref{eq:d_i}. Indeed, 
$$
cc\mu  \ast c\psi t \ast \psi ct \ast t\delta t \ast t \psi \ast t\eta c
= cc\mu  \ast \delta tt \ast \psi t \ast t \psi \ast t\eta c 
= \delta t \ast \psi.
$$
The first equality follows by \eqref{eq:wentw_d} and the second one follows by
\eqref{eq:wentw_m} and unitality of the monad $t$. 

Hence if moreover $\cK$ admits EM constructions for monads and idempotent
2-cells in $\cK$ split (hence there exists the pseudo-functor $J^w$) then, by
Corollary \ref{cor:comonad_lift}(1), the comonad $c$ has a weak
$\iota$-lifting for the monad $t$.  

For a commutative, associative and unital ring $k$, consider a $k$-algebra $A$
and a $k$-coalgebra $C$. Let $\Psi: C \otimes_k A \to A \otimes_k C$ be a
$k$-module map such that the triple $((-) \otimes_k A, (-)\otimes_k C,
(-)\otimes_k \Psi)$ is a weak entwining structure in $\mathrm{CAT}$ .
(If we are ready to cope with the more involved situation of a bicategory, we
can say simply that $(A,C,\Psi)$ is a weak entwining structure in
$\mathrm{BIM}$.) 
The corresponding weak $\iota$-lifting of the comonad $(-)\otimes_k C$ for the
monad $(-)\otimes_k A$ is studied in \cite[Section 2]{CaeDeG}. 
Brzezi\'nski showed in \cite[Proposition 2.3]{Brz:coring} that it can be
described as a comonad $(-)\otimes_A {\overline C}$ on the category of right
$A$-modules, where  
the $A$-coring (i.e. comonad $A \1c{} A$ in $\mathrm{BIM}$)
${\overline C}$ is constructed 
as a $k$-module retract of $A \otimes_k C$.  

Examples of weak entwining structures, thus examples of weak $\iota$-liftings
of comonads for monads, are provided by weak Doi-Koppinen data in \cite{Bohm}
(see \cite{CaeDeG}), i.e. by comodule algebras and module coalgebras of weak
bialgebras.  
Further examples are weak comodule algebras of bialgebras in \cite[Proposition
  2.3]{CaeJan:CommA}.    
\end{example}

\begin{example}\label{ex:partial_entw}
Another generalization of a mixed distributive law, motivated by partial
coactions of Hopf algebras, is due to Caenepeel and Janssen.
Following \cite[Proposition 2.6]{CaeJan:web}, a {\em partial entwining
  structure} in a 2-category $\cK$ consists of a monad $(k \1c t k,\mu ,\eta
)$, a comonad $(k\1c c k,\delta ,\varepsilon )$ and a 2-cell $tc \2c \psi ct$
in $\cK$, such that identities \eqref{eq:e_p} and \eqref{eq:SK_9} hold,
together with 
\begin{equation}
cc\mu  \ast c\psi t \ast c\eta ct \ast \delta t \ast \psi 
= c\psi \ast \psi c \ast t\delta .
\label{eq:SK_12}\\
\end{equation}
Observe that axiom \eqref{eq:SK_12} implies \eqref{eq:d_p}:
\begin{eqnarray*}
cc\mu  \ast c\psi t \ast \psi ct \ast \eta cct \ast \delta t \ast \psi
&=& cc\mu  \ast cc\mu t \ast c\psi tt \ast c\eta c tt \ast \delta tt \ast \psi
t \ast \eta ct \ast \psi\\ 
&=& cc\mu  \ast c \psi t \ast c\eta ct \ast \delta t \ast \psi 
= c\psi \ast \psi c \ast t\delta .  
\end{eqnarray*}
The first and last equalities follow by \eqref{eq:SK_12} and 
the second equality is obtained using associativity of $\mu $ and
\eqref{eq:l_proj}. This implies that 
$((c,\psi), c\psi \ast \psi c \ast \eta cc \ast \delta , \eta  \ast
\varepsilon )$ is a comonad in $\EM^w(\cK)$.
Thus if moreover $\cK$ is a 2-category which admits EM constructions for monads
and in which idempotent 2-cells split, then
we conclude by Corollary \ref{cor:comonad_lift}(2) that a partial entwining
structure $(t,c,\psi)$ in $\cK$ induces a weak $\pi$-lifting of the comonad $c$
for the monad $t$. 

Consider the particular case when
a monad $t:=(-)\otimes_k A$ in $\mathrm{CAT}$ is induced by an algebra $A$
over a commutative, associative and unital ring $k$, a comonad
$c:=(-)\otimes_k C$ is induced by a $k$-coalgebra $C$ and a natural
transformation $tc \2c \psi ct$ is induced by a $k$-module map $C \otimes_k A
\to A \otimes_k C$. Then the weak $\pi$-lifting of the comonad $c$ for the monad
$t$, induced by a partial entwining $\psi$, is a comonad 
$(-)\otimes_A {\overline C}$ on the category of right $A$-modules. The
$A$-coring ${\overline C}$ was constructed in \cite[Proposition
  2.6]{CaeJan:web} as a $k$-module retract of $A\otimes_k C$.

Examples of partial entwining structures (hence of weak $\pi$-liftings of a
comonad for a monad) are provided by partial comodule algebras of bialgebras in
\cite[Proposition 2.6]{CaeJan:CommA}. 
\end{example}

\begin{example}\label{ex:lax_entw}
Yet another way to generalize a mixed distributive law was proposed in
\cite{CaeJan:web}. Following \cite[Proposition 2.5]{CaeJan:web}, a {\em lax
entwining structure} in a 2-category $\cK$ consists of a monad $(k \1c t
k,\mu ,\eta )$, a comonad $(k\1c c k,\delta ,\varepsilon )$ and a 2-cell $tc
\2c \psi ct$ in $\cK$, such that identities \eqref{eq:SK_9},
\eqref{eq:wentw_e} and \eqref{eq:SK_12} hold, together with  
$$
c\mu  \ast ct\varepsilon t \ast ct \psi \ast ct\eta c \ast \psi c \ast t\delta
\ast \eta c = \psi \ast \eta c. 
$$
As we observed in Example \ref{ex:partial_entw}, \eqref{eq:SK_12} implies
\eqref{eq:d_p}, and \eqref{eq:wentw_e} is identical to
\eqref{eq:e_i}. However, none of
\eqref{eq:d_i} and \eqref{eq:e_p} seems to hold for an arbitrary lax entwining 
structure.  
Still, 
the axioms of a lax entwining structure allow us to prove that there is
a comonad 
$((c,\psi), c\psi \ast \psi c \ast \eta cc \ast \delta , \varepsilon
t \ast \psi \ast \eta c)$ in $\EM^w(\cK)$.
Therefore, if $\cK$ admits EM constructions for monads and idempotent 2-cells
in $\cK$ split, then $J^w$ takes it to a comonad
$(J^w(t) \1c {\overline c} J^w(t), \plift \delta , \ilift 
\varepsilon )$ in $\cK$. However, it is neither a weak $\iota $-lifting nor a
weak $\pi$-lifting of the comonad $c$, it is of a mixed nature.

In the particular case when a lax entwining structure in $\mathrm{CAT}$ is
induced by bimodules, the comonad $({\overline c}, \plift \delta , \ilift
\varepsilon )$ is induced by an $A$-coring, which was computed in
\cite[Proposition  2.5]{CaeJan:web}. 
Examples of lax entwining structures are provided by lax comodule algebras
of bialgebras in \cite[Proposition 2.5]{CaeJan:CommA}. 

A fourth logical possibility, to obtain a comonad structure on a weak lifting
for a monad $t$ of a 1-cell $c$ underlying a comonad $(c,\delta ,\varepsilon
)$, is to allow the comultiplication to be a weak $\iota$-lifting of $\delta $
and the counit to be a weak $\pi$-lifting of $\varepsilon$. That is, to
require a 1-cell $t \1c {(c,\psi)} t$ in $\EM^w(\cK)$ to satisfy
\eqref{eq:d_i} and \eqref{eq:e_p}. By (the proof of) Lemma
\ref{lem:lift_comp}, it yields a coassociative and counital comonad
$({\overline c},\ilift \delta , \plift \varepsilon )$ in $\cK$.
\end{example}

For any 2-category $\cK$, one may consider the vertically-opposite 2-category
$\cK_*$. The 2-category $\cK_*$ has the same 0-, 1-, and 2-cells as $\cK$, the
same horizontal composition and the opposite vertical composition. Obviously,
2-cells in $\cK$ split if and only if 2-cells in $\cK_*$ split. 
Since monads in $\cK_*$ are the same as the comonads in $\cK$, the 2-category
$\cK_*$ admits EM constructions for monads if and only if $\cK$ admits
EM constructions for comonads, cf. \cite{PowWat}. 
In this case we denote by $J^w_*: \EM^w(\cK_*)_* \to \cK$ the pseudo-functor in
Theorem \ref{thm:J^w}.

\begin{definition}
Consider a 2-category $\cK$ which admits EM constructions for comonads. 
We say that a 1-cell $V$ in $\cK$ possesses a {\em weak lifting ${\overline
    V}$ for some comonads} $c$ and $c'$, provided that, regarded as
1-cells in $\cK_*$, ${\overline V}$ is a weak lifting of $V$ for the monads
$c$ and $c'$ in $\cK_*$.

For a 2-cell $\omega$ in $\cK$, a {\em weak $\iota$-lifting} (resp. {\em weak
  $\pi$-lifting}) for some comonads $c$ and $c'$ in $\cK$ is defined
as a weak $\pi$-lifting (resp. weak $\iota$-lifting) of $\omega$, regarded as
a 2-cell in $\cK_*$, for the monads $c$ and $c'$ in $\cK_*$. 
\end{definition}

The following corollary is obtained by applying Corollary
\ref{cor:comonad_lift} to the vertically-opposite of a 2-category. 
Therein, the symbol $\ast$ denotes the vertical composition in $\cK$ (not its
opposite). 

\begin{corollary}\label{cor:monad_lift}
Consider a 2-category $\cK$ which admits EM constructions for comonads and in
which idempotent 2-cells split.  
Let $(k\1c t k, \mu ,\eta )$ be a monad and $(k\1c c k,\delta ,\varepsilon )$
be a comonad in $\cK$.  
\begin{itemize}
\item[{(1)}] The following assertions are equivalent.
\begin{itemize}
\item[{(i)}] There is a monad $((t,\psi),\mu ,\eta )$ in $\Mnd^\iota
  (\cK_*)_*$.  
That is, there exists a 1-cell $c \1c {(t,\psi)} c$ in $\EM^w(\cK_*)_*$
(i.e. a 2-cell $tc \2c \psi ct$ in $\cK$ such that 
$\delta t \ast \psi=c\psi \ast \psi c \ast t\delta $), satisfying 
\begin{eqnarray}
&& \psi \ast \mu c
= c\varepsilon t \ast c\psi \ast c\mu c \ast \psi tc \ast t \psi c \ast
tt\delta ; 
\label{eq:m_p}\\
&& \psi \ast \eta c 
= c\varepsilon t \ast c\psi \ast c\eta c \ast \delta .
\label{eq:u_p}
\end{eqnarray}
 \item[{(ii)}] There is a monad
$(c \1c {(t,\psi)} c, \varepsilon t \ast \psi \ast \mu c, \varepsilon t \ast
   \psi \ast \eta c)$ in $\EM^w(\cK_*)_*$.
\item[{(iii)}] There is a monad $(J_*^w(c) \1c {\overline{t}} J^w_*(c), \plift
  \mu , \plift \eta )$ in $\cK$ such that
  $\plift \mu $ is a weak $\pi$-lifting of $\mu $ and $\plift \eta $ is a weak
  $\pi$-lifting of $\eta $. 
\end{itemize}
If these equivalent statements hold, then we say shortly that the monad
$(\overline{t}, \plift \mu , \plift \eta )$ is a weak $\pi$-lifting of the
monad $(t,\mu ,\eta )$ for the comonad $(c,\delta ,\varepsilon )$.

\item[{(2)}] The following assertions are equivalent.
\begin{itemize}
\item[{(i)}] There is a monad $((t,\psi),\mu ,\eta )$ in $\Mnd^\pi (\cK_*)_*$. 
That is, there exists a 1-cell $c \1c {(t,\psi)} c$ in $\EM^w(\cK_*)_*$,
  satisfying 
\begin{eqnarray}
&& c\mu  \ast \psi t \ast t \psi 
= \psi \ast \mu c \ast \varepsilon ttc \ast \psi tc \ast t \psi c \ast
tt\delta ; 
\label{eq:m_i}\\
&& c\eta  = \psi \ast \eta c. 
\label{eq:u_i}
\end{eqnarray}
 \item[{(ii)}] There is a monad
$(c \1c {(t,\psi)} c, \mu  \ast \varepsilon tt \ast \psi t \ast t \psi, \eta
   \ast \varepsilon )$ in $\EM^w(\cK_*)_*$.
\item[{(iii)}] There is a monad $(J^w_*(c) \1c {\overline{t}} J^w_*(c), \ilift
  \mu , \ilift \eta )$ in $\cK$ such that
  $\ilift \mu $ is a weak $\iota$-lifting of $\mu $ and $\ilift \eta $ is a
  weak $\iota$-lifting of $\eta $. 
\end{itemize}
If these equivalent statements hold, then we say shortly that the monad
$(\overline{t}, \ilift \mu , \ilift \eta )$ is a weak $\iota$-lifting of the
monad $(t,\mu ,\eta )$ for the comonad $(c,\delta ,\varepsilon )$.

\item[{(3)}] The following assertions are equivalent.
\begin{itemize}
\item[{(i)}] There exists a 1-cell $t \1c {(c,\psi)} t$ in $\EM^w(\cK_*)_*$,
  satisfying 
\begin{eqnarray}
&& c\mu  \ast \psi t \ast t \psi 
= \psi \ast \mu c;
\label{eq:m_b}\\
&& c\eta  = \psi \ast \eta c. 
\label{eq:u_b}
\end{eqnarray}
\item[{(ii)}] There is a monad $(J^w_*(c) \1c {\overline t} J^w_*(c), \overline
  \mu , \overline \eta )$ in $\cK$ such that
  $\overline \mu $ is both a weak $\iota$-lifting and a weak $\pi$-lifting of
  $\mu $ and $\overline \eta $ is both a weak $\iota$-lifting and a weak
  $\pi$-lifting of $\eta $.   
\end{itemize}   
\end{itemize}
\end{corollary}

A 2-cell $\psi$ in Corollary \ref{cor:monad_lift}(3)(i) differs from a mixed
distributive law by the compatibility condition with the counit of the
comonad. 

 In a 2-category $\cK$ which admits EM constructions for both monads and
 comonads and in which idempotent 2-cells split, one can say more about weak
 entwining structures than it was said in Example \ref{ex:weak_entw}.

\begin{proposition} \label{prop:weak_entw}
Consider a 2-category $\cK$ which admits EM constructions for both monads and
comonads and in which idempotent 2-cells split. 
For a monad $(k \1c t k,\mu ,\eta )$, a comonad $(k\1c c k, \delta ,
\varepsilon )$, and a 2-cell $tc \2c \psi ct$ in $\cK$, the following
assertions are equivalent. 
\begin{itemize}
\item[{(i)}] The triple $(t,c,\psi)$ is a weak entwining structure,
  i.e. axioms \eqref{eq:wentw_m} -- \eqref{eq:wentw_e} are satisfied.
\item[{(ii)}] The 2-cell $\psi$ induces both a weak $\iota$-lifting of the
  comonad $c$ for the monad $t$ and a weak $\pi$-lifting of the monad $t$ for
  the comonad $c$. That is to say, the assertions in Corollary
  \ref{cor:comonad_lift}(1) and Corollary \ref{cor:monad_lift}(1) hold. 
\end{itemize}
\end{proposition}

\begin{proof}
We have seen in Example \ref{ex:weak_entw} that axioms
\eqref{eq:wentw_m}-\eqref{eq:wentw_e} imply \eqref{eq:d_i}. 
Similarly, they can be seen to imply \eqref{eq:m_p} as well, applying first
\eqref{eq:wentw_m} and next \eqref{eq:wentw_d}.
\end{proof}

Proposition \ref{prop:weak_entw} is the basis of a construction in
\cite{wentw} of a
2-category of weak entwining structures in any 2-category. In that paper, for
a weak entwining structure in a 2-category $\cK$ which admits EM constructions
for both monads and comonads and in which idempotent 2-cells split, it is
proven that the weakly lifted monad, and the weakly lifted comonad, occurring
in part (ii) of Proposition \ref{prop:weak_entw}, possess equivalent
Eilenberg-Moore objects. 

The characterization of weak entwining structures in Proposition
\ref{prop:weak_entw} can be used, in particular,
do describe weak bialgebras \cite{BNSz} in terms of weak liftings. Recall that
a {\em weak bialgebra} $H$ over a commutative, associative and unital ring
$k$, is a $k$-module which possesses both a $k$-algebra and a $k$-coalgebra
structure, subject to the following compatibility conditions.  
Denote the multiplication $H \otimes_k H \to H$ in $H$ by juxtaposition of
elements. 
Write $1$ for the unit element of the algebra $H$ and write $\varepsilon:H
\to k$ for the counit.
For the comultiplication $H \to H \otimes_k H$, use a Sweedler type index
notation $h \mapsto \sum h_1 \otimes_k h_2$.
With these notations, the axioms 
\begin{eqnarray}
&&\sum (hh')_1 \otimes_k (hh')_2 
= \sum h_1 h'_1 \otimes_k h_2 h'_2;
\label{eq:WBA_1}\\
&&\sum 1_1 \otimes 1_2 1_{1'} \otimes 1_{2'}
= \sum 1_1 \otimes_k 1_2 \otimes_k 1_3
= \sum 1_1 \otimes_k 1_{1'} 1_2 \otimes_k 1_{2'};
\label{eq:WBA_2}\\
&& \sum \varepsilon(h1_1) \varepsilon(1_2 h')
=\varepsilon(hh')
= \sum \varepsilon(h 1_2) \varepsilon(1_1h')
\label{eq:WBA_3}
\end{eqnarray}
are required to hold, for all elements $h$ and $h'$ of $H$.
However, the comultiplication is not required to preserve the unit, i.e. $\sum
1_1 \otimes_k 1_2$ is not required to be equal to $1 \otimes 1$ and the counit
is not required to be multiplicative, i.e. $\varepsilon(hh')$ is not required to
be equal to $\varepsilon(h) \varepsilon(h')$, for elements $h,h'\in H$.

In the following proposition we deal with the (co)monads $H \otimes_k (-)$ and
$(-)\otimes_k H$, induced by a $k$-(co)algebra $H$ on the category of
modules over a commutative, associative and unital ring $k$. 

\begin{proposition}
For a commutative, associative and unital ring $k$, consider a $k$-module $H$
which possesses both a $k$-algebra and a $k$-coalgebra structure. 
Using the notations introduced above the proposition, the following assertions
are equivalent. 
\item[{(i)}] The algebra and coalgebra structures of $H$ constitute a weak
  bialgebra; 
\item[{(ii)}] The $k$-module map
\begin{equation}\label{eq:right}
\Psi_R:H \otimes_k H \to H\otimes_k H, \qquad
h \otimes_k h' \mapsto \sum h'_1 \otimes_k hh'_2
\end{equation}
induces a weak $\iota$-lifting of the comonad $(-) \otimes_k H$ for the monad
$(-)\otimes_k H$ and a weak $\pi$-lifting of the monad $(-)\otimes_k H$ for the
comonad $(-)\otimes_k H$, and the $k$-module map 
\begin{equation}\label{eq:left}
\Psi_L:H \otimes_k H \to H\otimes_k H, \qquad 
h \otimes_k h' \mapsto \sum h_1 h' \otimes_k h_2
\end{equation}
induces a weak $\iota$-lifting of the comonad $H \otimes_k (-)$ for the monad
$H \otimes_k (-)$ and a weak $\pi$-lifting of the monad $H\otimes_k (-)$ for
the comonad $H \otimes_k (-)$. That is to say,
$$
((-)\otimes_k H, (-)\otimes_k \Psi_R)\
\qquad \textrm{and}\qquad 
(H \otimes_k (-), \Psi_L\otimes_k(-))
$$
are comonads in $\Mnd^\iota (\mathrm{CAT})$, via the comultiplication and counit
induced by the coalgebra $H$, and they
are monads in $\Mnd^\iota (\mathrm{CAT}_*)_*$, via the multiplication and unit
induced by the algebra $H$. 
\end{proposition}

\begin{proof}
Note first that assertion (ii) implies \eqref{eq:WBA_1}. Indeed,
\eqref{eq:right} determines a 1-cell $((-)\otimes_k H, (-) \otimes_k \Psi_R)$
in $\EM^w(\mathrm{CAT})$ if and only if
$$
\sum (h'h'')_1 \otimes_k h (h'h'')_2 
= \sum h'_1 h''_1 \otimes_k h h'_2 h''_2,
$$
for any elements $h,h'$ and $h''$ of $H$. Putting $h=1$ we obtain
\eqref{eq:WBA_1}. 

By Proposition \ref{prop:weak_entw}, assertion (ii) is equivalent to saying
that 
$((-) \otimes_k H, (-)\otimes_k H, (-)\otimes_k \Psi_R)$ and 
$(H \otimes_k (-), H \otimes_k (-), \Psi_L \otimes_k(-))$ are weak entwining
structures in $\mathrm{CAT}$ (or, in the terminology of \cite{CaeDeG},
$(H,H,\Psi_R)$ is a right-right weak entwining structure and $(H,H,\Psi_L)$ is
a left-left weak entwining structure in $\mathrm{BIM}$).
This statement was proven to be equivalent to (i) in \cite[Theorem
  4.7]{CaeDeG}. 
\end{proof}

\section*{Acknowledgement}
The author's work is supported by the Hungarian Scientific Research Fund OTKA
F67910. She would like to thank the referee for valuable and
  helpful comments.


\begin{thebibliography}{Bibliography}{}

\bibitem{Beck} J. Beck, {\em Distributive laws}, [in:] {\em Seminar on Triples 
 and Categorical Homology Theory}, 
 B. Eckmann (ed.), Springer LNM 80, 119-140 (1969).    

\bibitem{BCM} R. Blattner, M. Cohen and S. Montgomery, {\em Crossed products
  and inner actions of Hopf algebras}, Trans. of Amer. Math. Soc. 298 (1986),
  671-711. 

\bibitem{Bohm} G. B\"ohm, {\em Doi-Hopf Modules over Weak Hopf Algebras},
  Comm. Algebra 28 (2000), 4687-4698.

\bibitem{wentw} G. B\"ohm, {\em The 2-category of weak entwining structures},
\href{http://arxiv.org/abs/0902.4197}
{arXiv:0902.4197v1}

\bibitem{BB} G. B\"ohm and T. Brzezi\'nski, {\em Cleft extensions of Hopf
  algebroids}, Appl. Categ. Structures 14 (2006), 411-419. 
{\em Corrigendum}, Appl. Categ. Structures 17 (2009), 613-620.

\bibitem{BNSz} G. B\"ohm, F. Nill and K. Szlach\'anyi, {\em Weak Hopf
  algebras I: Integral theory and $C^*$-structure}, J. Algebra {221}
  (1999), 385--438.
 
\bibitem{Brz:crosp} T. Brzezi\'nski, {\em Crossed products by a coalgebra},
  Comm. in Algebra 25 (1997), 3551-3575. 

\bibitem{Brz:coring}  T. Brzezi\'nski, {\em The structure of
  corings. Induction functors, Maschke-type theorem, and Frobenius and
  Galois-type properties},  Algebr. and Representat. Theory 5 (2002),
  389-410. 

\bibitem{CaeDeG} S. Caenepeel and E. De Groot, {\em Modules over weak
  entwining structures}, Contemp. Math. 267 (2000), 4701--4735.

\bibitem{CaeJan:web} S. Caenepeel and K. Janssen, {\em Partial entwining
  structures}, available at 
\href{http://homepages.vub.ac.be/~scaenepe/quotient10.pdf}
{\tt http://homepages.vub.ac.be/~scaenepe/quotient10.pdf}.   

\bibitem{CaeJan:CommA} S. Caenepeel and K. Janssen, {\em Partial (co)actions
  of Hopf algebras and partial Hopf-Galois theory}, Comm. in Algebra 36 (2008),
  2923-2946.  

\bibitem{DT} Y. Doi and M. Takeuchi, {\em Cleft comodule algebras for a
  bialgebra}, Comm. in Algebra 14 (1986), 801-817.

\bibitem{WCP} J.M. Fern\'andez Vilaboa, R. Gonz\'alez Rodr\'{\i}guez and
  A.B. Rodr\'{\i}guez Raposo, {\em Preunits and weak crossed products},
J. Pure and Applied Algebra 213 (2009), 2244-2261.

\bibitem{Pepe} J. G\'omez-Torrecillas, {\em Comonads and Galois Corings},
Applied Categorical Structures 14 (2006), 579-598.

\bibitem{KeSt} G.M. Kelly and R. Street,
  {\em Review of the elements of 2-categories}, [in:]
  {\em Category Sem. Proc. Sydney 1972/1973}, G.M. Kelly (ed.), Springer
  LNM 420, 75-103 (1974).  

\bibitem{LS} S. Lack and R. Street, {\em The formal theory of monads II},
  J. Pure and Applied Algebra 175 (2002), 243-265.

\bibitem{PowWat} J. Power and H. Watanabe, {\em Combining a monad and a
  comonad}, Theoret. Comput. Sci. 280 (2002), 137-162.

\bibitem{RodRap} A.B. Rodr\'{\i}guez Raposo, {\em Crossed products for weak
  Hopf algebras}, Comm. in Algebra 37 (2009), 2274-2289.

\bibitem{Street} R. Street, {\em The formal theory of monads}, J. Pure and
  Applied Algebra 2 (1972), 149-168.

\bibitem{priv} R. Street, {\em Private communication}, 2008.

\bibitem{Swe} M.E. Sweedler, {\em Cohomology of algebras over Hopf algebras},
  Trans. Amer. Math. Soc. 133 (1968), 205-239. 

\bibitem{Wis} R. Wisbauer, {\em Algebras Versus Coalgebras}, Applied
  Categorical Structures 16 (2008), 255-295.

\end{thebibliography}
\end{document}